\documentclass[12pt]{amsart}
\usepackage{amsmath,fullpage}
\usepackage{amscd}
\usepackage{amssymb}

\usepackage[all,cmtip]{xy}

\newtheorem{theorem}{Theorem}[section]
\newtheorem{corollary}[theorem]{Corollary}
\newtheorem{lemma}[theorem]{Lemma}
\newtheorem{proposition}[theorem]{Proposition}

\theoremstyle{definition}

\newtheorem{remark}[theorem]{Remark}

\newcommand{\tk}{{\textbf{k}}}
\newcommand{\cB}{{\mathcal B}}

\newcommand{\cO}{{\mathcal O}}

\newcommand{\cF}{{\mathcal F}}
\newcommand{\Hom}{{\text{Hom}}}

\newcommand{\End}{{\text{End}}}
\newcommand{\Ad}{{\text{Ad}}}
\newcommand{\bQ}{{\mathbb Q}}

\newcommand{\rc}{{\mathrm c}}
\newcommand{\Lb}{{\mathfrak{b}}}
\newcommand{\Lt}{{\mathfrak{t}}}
\newcommand{\tF}{{\textbf{F}}}

\newcommand{\Lg}{{\mathfrak g}}
\newcommand{\Lo}{{\mathfrak o}}
\newcommand{\A}{{\textbf{A}}}

\newcommand{\supp}{{\text{supp}}}
\newcommand{\diag}{{\text{diag}}}

\newcommand{\cH}{{\mathcal H}}
\newcommand{\cT}{{\mathcal T}}

\newcommand{\Ln}{{\mathfrak{n}}}

\begin{document}
\title[Nilpotent orbits in classical Lie algebras ]{Nilpotent orbits
in classical Lie algebras over finite fields of characteristic 2 and
the Springer correspondence}
        \author{Ting Xue}
        \address{Department of Mathematics, Massachusetts Institute of Technology,
Cambridge, MA 02139, USA}
        \email{txue@math.mit.edu}
   
   \maketitle     
\begin{abstract}
Let $G$ be an adjoint algebraic group of type $B$, $C$, or $D$ over an algebraically closed field of characteristic 2. We construct a Springer correspondence for the Lie algebra of $G$. In particular, for orthogonal Lie algebras in characteristic 2, the structure of component groups of nilpotent centralizers is determined and the number of nilpotent orbits over finite fields is obtained.
\end{abstract}

\section{ Introduction}

Throughout this paper, $\tk$ denotes an algebraically closed field
of characteristic 2, $\tF_q$ denotes a finite field of
characteristic 2 and $\bar{\tF}_q$ denotes an algebraic closure of
$\tF_q$.

In \cite{Hes}, Hesselink determines the nilpotent orbits in
classical Lie algebras under the adjoint action of classical Lie
groups over $\tk$. In \cite{Spal}, Spaltenstein gives a
parametrization of these nilpotent orbits by pairs of partitions. We
extend Hesselink's method to study the nilpotent orbits in the Lie
algebras of orthogonal groups over $\tF_q$. Using this extension and
Spaltenstein's parametrization we classify the nilpotent orbits over
$\tF_q$. We determine the structure of the component groups of
centralizers of nilpotent elements. In particular, we obtain the
number of nilpotent orbits over $\tF_q$.

Let $G$ be a connected reductive algebraic group over an
algebraically closed field and $\Lg$ the Lie algebra of $G$. When
the characteristic of the field is large enough, Springer \cite{Sp1}
constructs representations of the Weyl group of $G$ which are
related to the nilpotent $G$-orbits in $\Lg$. Lusztig \cite{Lu1}
constructs the generalized Springer correspondence in all
characteristic which is related to the unipotent conjugacy classes
in $G$. Let $G_{ad}$ be an adjoint algebraic group of type $B,C$ or
$D$ over $\tk$ and $\Lg_{ad}$ be the Lie algebra of $G_{ad}$.
 We use a similar construction as in
\cite{Lu1,Lu4} to give the Springer correspondence for $\Lg_{ad}$.
Let $\mathfrak{A}_{ad}$ be the set of all pairs
$(\mathrm{c},\mathcal{F})$ where $\mathrm{c}$ is a nilpotent
$G_{ad}$-orbit in $\Lg_{ad}$ and $\mathcal{F}$ is an irreducible
$G_{ad}$-equivariant local system on $\mathrm{c}$ (up to
isomorphism). We construct a bijective map from the set of
isomorphism classes of irreducible representations of the Weyl group
of $G_{ad}$ to the set $\mathfrak{A}_{ad}$. In the case of
symplectic group a Springer correspondence (with a different
definition than ours) has been established in \cite{kato}; in that
case centralizers of the
 nilpotent elements are connected \cite{Spal}.
A complicating feature in the orthogonal case is the existence of
non-trivial equivariant local systems on a nilpotent orbit.

Most of the results in this paper have been announced in \cite{X}.

\section{ Hesselink's classification of nilpotent orbits over an
algebraically closed field}

We recall the results of Hesselink about nilpotent orbits in
orthogonal algebras in this section (see \cite{Hes}). Let
$\mathbb{K}$ be a field of characteristic 2, not necessarily
algebraically closed.
\subsection{}

 A form space $V$ is a finite
dimensional vector space over $\mathbb{K}$ equipped with a quadratic
form $Q: V\rightarrow \mathbb{K}$. Let $\langle\cdot,\cdot\rangle:
V\times V\rightarrow \mathbb{K}$ be the bilinear form $\langle
v,w\rangle=Q(v+w)+Q(v)+Q(w)$. Let $V^\perp=\{v\in V|\langle
v,w\rangle=0,\forall\ w\in V\}$. A form space $V$ is called
non-defective if $V^\perp=\{0\}$, otherwise, it is called defective.
$V$ is called non-degenerate if
$\text{dim}(V^\perp)\leq 1$ and $Q(v)\neq 0$ for all non-zero $v\in
V^{\perp}$.

Let $V$ be a non-degenerate form space of dimension $N$ over
$\mathbb{K}$. Define the orthogonal group $O(V)$ to be
$\{g\in\text{GL}(V)|Q(gv)=Q(v),\ \forall\ v\in V\}$ and define
$\mathfrak{o}(V)$ to be $\{x\in\End (V)|\langle xv,v\rangle=0,\
\forall\ v\in V\text{ and }\text{tr}(x)=0\}$.  In the case where $\mathbb{K}$ is algebraically
closed, let $SO(V)$ be the identity component of $O(V)$ and we define
$O_{N}(\mathbb{K})$, $\Lo_N(\mathbb{K})$ and $SO_{N}(\mathbb{K})$ to be $O(V)$, $\Lo(V)$  and $SO(V)$ respectively.

\subsection{}\label{ssec-1} A form module is defined to be a pair $(V,T)$ where $V$
is a non-degenerate form space and $T$ is a nilpotent element in
$\mathfrak{o}(V)$. Classifying nilpotent orbits in $\mathfrak{o}(V)$
is equivalent to classifying form modules $(V,T)$ on the form space
$V$. Let $A=\mathbb{K}[[t]]$ be the ring of formal power series in
the indeterminate $t$. The form module $V=(V,T)$ is considered as an
$A$-module by $(\sum_{n\geq 0} a_nt^n)v=\sum_{n\geq 0} a_nT^nv$.

Let $E$ be the vector space spanned by the linear functionals
$t^{-n}: A\rightarrow\mathbb{K},\ \sum a_it^i\mapsto a_n,\ n\geq 0$.
Let $E_0$ be the subspace $\sum_{n\geq 0} \mathbb{K} t^{-2n}$ and
$\pi_0:E\rightarrow E_0$ be the natural projection. The space $E$ is
considered as an $A$-module by $(au)(b)=u(ab)$, for $a,b\in A,u\in
E$.

An abstract form module is defined to be an $A$-module $V$ with
$\text{dim}(V)<\infty$, which is equipped with mappings
$\varphi:V\times V\rightarrow E$ and $\psi: V\rightarrow E_0$
satisfying the following axioms:
\begin{enumerate}
  \item[(i)] The map $\varphi(\cdot,w)$ is $A$-linear for every $w\in V$.
  \item[(ii)] $\varphi(v,w)=\varphi(w,v)$ for all $v,w\in V.$
  \item[(iii)] $\varphi(v,v)=0$ for all $v\in V$.
  \item[(vi)] $\psi(v+w)=\psi(v)+\psi(w)+\pi_0(\varphi(v,w))$ for all $v,w\in V.$
  \item[(v)] $\psi(a v)=a^2\psi(v)$ for $v\in V,\ a\in A$.
\end{enumerate}
The following proposition identifies a form module $(V,T)$ and the
corresponding abstract form module $V=(V,\varphi,\psi)$.
\begin{proposition}(\cite{Hes})
If $(V,\varphi,\psi)$ is an abstract form module, then
$(V,\langle\cdot,\cdot\rangle,Q)$ given by $\mathrm{(i)}$ is a form
module. If $(V,\langle\cdot,\cdot\rangle,Q)$ is a form module, there
is a unique abstract form module $(V,\varphi,\psi)$ such that
$\mathrm{(i)}$ holds; it is given by $\mathrm{(ii)}$.
\begin{eqnarray*}
&&\mathrm{(i)}\ \langle v,w\rangle=\varphi(v,w)(1),\
Q(v)=\psi(v)(1).\\
&&\mathrm{(ii)}\ \varphi(v,w)=\sum_{n\geq 0}\langle t^nv,w\rangle
t^{-n},\ \psi(v)=\sum_{n\geq 0} Q(t^nv)t^{-2n}.
\end{eqnarray*}
\end{proposition}

\subsection{} An element in $\mathfrak{o}(V)$ is nilpotent if and
only if it is nilpotent in $\End (V)$. Let $T$ be a nilpotent
element in $\Lo(V)$. There exists a unique sequence of integers
$p_1\geq\cdots\geq p_s\geq 1$ and a family of vectors
$v_1,\ldots,v_s$ such that $T^{p_i}v_i=0$ and the vectors
$T^{q_i}v_i$, $0\leq q_i\leq p_i-1$ form a basis of $V$. We write
$p(V,T)=(p_1,\ldots,p_s)$. Define the index function $\chi(V,T):
\mathbb{N}\rightarrow\mathbb{Z}$ by $\chi(V,T)(m)=\min\{k\geq
0|T^mv=0\Rightarrow Q(T^kv)=0\}.$ Define $\mu(V)$ to be the minimal
integer $m\geq 0$ such that $t^m V=0$. For $v\in V$ (or $E$), we
define $\mu(v)=\min\{m\geq 0|t^m v=0\}$.

\subsection{}\label{sec-1} Let $V$ be a form module. An orthogonal
decomposition of $V$ is an expression of $V$ as a direct sum
$V=\sum_{i=1}^r V_i$ of mutually orthogonal submodules $V_i$. The
form module $V$ is called indecomposable if $V\neq 0$ and for every
orthogonal decomposition $V=V_1\oplus V_2$ we have $V_1=0$ or
$V_2=0$. Every form module $V$ has some orthogonal decomposition
$V=\sum_{i=1}^r V_i$ in indecomposable submodules
$V_1,V_2,\ldots,V_r$. The indecomposable modules are classified as
follows.
\begin{proposition}\label{prop-indec}(\cite{Hes})
Let $V$ be a non-degenerate indecomposable form module. There exist
$v_1,v_2\in V$ such that $V=Av_1\oplus Av_2$ and
$\mu(v_1)\geq\mu(v_2)$. For any such pair we put
$m=\mu(v_1),m'=\mu(v_2),\Phi=\varphi(v_1,v_2)$ and
$\Psi_i=\psi(v_i)$. One of the following conditions holds:
\begin{enumerate}
  \item[(i)] $m'=\mu(\Phi)=m,\ \mu(\Psi_i)\leq 2m-1.$
  \item[(ii)] $m'=\mu(\Phi)=m-1,\ \mu(\Psi_1)=2m-1>\mu(\Psi_2).$
\end{enumerate}
Conversely, let $m\in\mathbb{N}$, $m'\in\mathbb{N}\cup\{0\}$,
${\Phi\in E}$, $\Psi_1,\Psi_2\in E_{0}$
 be given satisfying $\mathrm{(i)}$ or $\mathrm{(ii)}$.
 Up to a canonical isomorphism there exists a
unique form module $V=Av_1\oplus Av_2$ with
$m=\mu(v_1),m'=\mu(v_2),\Phi=\varphi(v_1,v_2)$ and
$\Psi_i=\psi(v_i)$. This form module is indecomposable. In case
$\mathrm{(i)}$ it is non-defective. In case $\mathrm{(ii)}$ it is
defective and non-degenerate.
\end{proposition}

From now on assume $\mathbb{K}$ is algebraically closed. The
indecomposable modules in Proposition \ref{prop-indec} are
normalized in \cite{Hes} section 3.4 and 3.5 as follows.
\begin{proposition}\label{prop-norm}(\cite{Hes})
The indecomposable non-degenerate form modules over $\mathbb{K}$ are

$\mathrm{(i)}$ $W_l(m)=Av_1\oplus Av_2,[\frac{m+1}{2}]\leq l\leq m$,
$\mu(v_1)=\mu(v_2)=m$, $\psi(v_1)=t^{2-2l}$, $\psi(v_2)=0$ and
$\varphi(v_1,v_2)=t^{1-m}$; ($[\alpha]$ means the integer part of
$\alpha$.)

$\mathrm{(ii)}$ $D(m)=Av_1\oplus Av_2$,  $\mu(v_1)=m,\mu(v_2)=m-1$,
$\psi(v_1)=t^{2-2m}$, $\psi(v_2)=0$ and $\varphi(v_1,v_2)=t^{2-m}$.

We have $\chi_{W_l(m)}=[m;l]$ and $\chi_{D(m)}=[m;m]$, where
$[m;l]:\mathbb{N}\rightarrow\mathbb{Z}$ is defined by
$[m;l](n)=\max\{0,\min\{n-m+l,l\}\}.$ Among these types of
indecomposable modules only the types $D(m)$ are defective.
\end{proposition}
\begin{remark}
The notation we use here is slightly different from that of
\cite{Hes}. The form module $W_{[\frac{m+1}{2}]}(m)$ in
$\mathrm{(i)}$ is isomorphic to the form module $W(m)$ in
\cite{Hes}.
\end{remark}
Finally this normalization of indecomposable modules is used to
classify all the form modules. Let $(V,T)$ be a non-degenerate from
module with
$p(V,T)=(\lambda_1,\ldots,\lambda_1,\ldots,\lambda_k,\ldots,\lambda_k)$
where $\lambda_1>\cdots>\lambda_k\geq 1$ and index function
$\chi=\chi(V,T)$. Let $m_i\in\mathbb{N}$ be the multiplicity of
$\lambda_i$ in $p(V,T)$. The isomorphism class of $(V,T)$ is
determined by the symbol
$$S(V,T)=(\lambda_1)_{\chi(\lambda_1)}^{m_1}(\lambda_2)_{\chi(\lambda_2)}^{m_2}
\cdots(\lambda_k)_{\chi(\lambda_k)}^{m_k}.$$ A symbol $S$ of the
above form is the symbol of an isomorphism class of non-degenerate
form modules if and only if it satisfies the following conditions

(i) $\chi(\lambda_i)\geq\chi(\lambda_{i+1})$ and
$\lambda_i-\chi(\lambda_i)\geq \lambda_{i+1}-\chi(\lambda_{i+1})$,
for $i=1,\ldots,k-1$;

(ii) $\frac{\lambda_i}{2}\leq\chi(\lambda_i)\leq\lambda_i$, for
$i=1,\ldots,k$;

(iii) $\chi(\lambda_i)=\lambda_i$ if $m_i$ is odd, for
$i=1,\ldots,k$;

(iv) $\{\lambda_i|m_i\text{ odd }\}=\{m,m-1\}\cap\mathbb{N}$ for
some $m\in\mathbb{Z}$.

In the following we denote by a symbol either a form module in the
isomorphism class or the corresponding nilpotent orbit.

\section{Indecomposable modules over $\tF_q$}

In this section, we study the non-degenerate indecomposable form
modules over $\tF_q$. Note that the classification of the
indecomposable modules (Proposition \ref{prop-indec}) is valid over
any field. Similar to \cite{Hes} Section 3.5, the non-degenerate
indecomposable form modules over $\tF_q$ are normalized as follows.
Fix an element $\delta$ in $\tF_q$ such that $\delta\notin
\{x^2+x|x\in \tF_q\}$.
\begin{proposition}\label{prop-nind}
The non-degenerate indecomposable form modules over $\tF_q$ are

$\mathrm{(i)}$ $W_l^0(m)=Av_1\oplus Av_2,\ [\frac{m+1}{2}]\leq l\leq
m$, with
 $\mu(v_1)=\mu(v_2)=m$, $\psi(v_1)=t^{2-2l},\psi(v_2)=0$ and $\varphi(v_1,v_2)=t^{1-m}$;

$\mathrm{(ii)}$ $W_l^\delta(m)=Av_1\oplus Av_2,\ \frac{m+1}{2}\leq
l\leq m$, with $\mu(v_1)=\mu(v_2)=m$, $\psi(v_1)=t^{2-2l}$,
$\psi(v_2)=\delta t^{2l-2m}$ and $\varphi(v_1,v_2)=t^{1-m}$;

$\mathrm{(iii)}$ $D(m)=Av_1\oplus Av_2$ with
$\mu(v_1)=m,\mu(v_2)=m-1$, $\psi(v_1)=t^{2-2m}$, $\psi(v_2)=0$ and
$\varphi(v_1,v_2)=t^{2-m}$.

We have $\chi_{W_l^0(m)}=\chi_{W_l^\delta(m)}=[m;l]$ and
$\chi_{D(m)}=[m;m]$. Among these types only the types $D(m)$ are
defective.
\end{proposition}
\begin{proof}
As pointed out in \cite{Hes}, the form modules in Proposition
\ref{prop-indec} (ii) over $\tF_q$ can be normalized the same as in
Proposition \ref{prop-norm} (ii). Namely, there exist $v_1$ and
$v_2$ such that the above modules have the form (iii).

Now let $U(m)=A v_1\oplus A{v}_2$ be a form module as in Proposition
\ref{prop-indec} (i) with $\mu(v_1)=\mu(v_2)=m$, we can assume
$\mu(\Psi_1)\geq\mu(\Psi_2)$. We have the following cases:

Case 1: $\Psi_1=\Psi_2=0$. We can assume $\Phi=t^{1-m}$. Let
$\tilde{v}_1=v_1+t^{m+1-2[\frac{m+1}{2}]}v_2$ and $\tilde{v}_2=v_2$.
We have $\psi(\tilde{v}_1)=t^{2-2[\frac{m+1}{2}]},\
\psi(\tilde{v}_2)=0,
         \ \varphi(\tilde{v}_1,\tilde{v}_2)=t^{1-m}$.

Case 2: $\Psi_1\neq 0,\ \Psi_2=0$. We can assume $\Phi=t^{1-m}$ and
$\Psi_1=t^{-2l}$, where $l\leq m-1$. If $l<[\frac{m-1}{2}]$, let
$\tilde{v}_1=v_1+t^{m+1-2[\frac{m+1}{2}]}v_2+t^{m-2l-1}v_2,\
\tilde{v}_2=v_2$; otherwise, let $\tilde{v}_1=v_1,\
\tilde{v}_2=v_2$. Then we get $\psi(\tilde{v}_1)=t^{-2l},\
[\frac{m-1}{2}]\leq l\leq m-1,\ \psi(\tilde{v}_2)=0,
         \ \varphi(\tilde{v}_1,\tilde{v}_2)=t^{1-m}$.

Case 3: $\Psi_1\neq 0,\ \Psi_2\neq 0$. We can assume
$\Psi_1=t^{-2l_1},\ \Phi=t^{1-m},\ \Psi_2=\sum_{i=0}^{l_2}a_i
t^{-2i}$ with $l_2\leq l_1\leq m-1$.

\begin{enumerate}
\item[(1)] $l_1<[\frac{m}{2}]$, let
    $\tilde{v}_2=v_2+\sum_{i=0}^{m-1}x_it^iv_1$, then
    $\psi(\tilde{v}_2)=0$ has a solution for $x_i$'s and we get Case 2.
    \item[(2)] $l_1\geq [\frac{m}{2}]$, let
    $\tilde{v}_2=v_2+\sum_{i=0}^{m-1}x_it^iv_1$.
     If $a_{m-l_1-1}\in\{x^2+x|x\in \tF_q\}$, then $\psi(\tilde{v}_2)=0$ has a solution for $x_i$'s and we get Case 2.
     If $a_{m-l_1-1}\notin\{x^2+x|x\in \tF_q\}$, then $\psi(\tilde{v}_2)=\delta t^{-2(m-l_1-1)}$ has a solution for $x_i$'s.
\end{enumerate}
Summarizing Cases 1-3, we have normalized $U(m)=Av_1\oplus Av_2$
with $\mu(v_1)=\mu(v_2)=m$ as follows:

\begin{enumerate}
             \item[(i)] $[\frac{m+1}{2}]\leq\chi(m)=l\leq m$, $\psi(v_1)=t^{2-2l},\ \psi(v_2)=0,\
             \varphi(v_1,v_2)=t^{1-m}$, denote by $W_l^0(m)$.
             \item[(ii)] $\frac{m+1}{2}\leq\chi(m)=l\leq m$, $\psi(v_1)=t^{2-2l},\ \psi(v_2)=\delta t^{-2(m-l)},\
\varphi(v_1,v_2)=t^{1-m}$, denote by $W_l^{\delta}(m)$.
\end{enumerate}
One can verify that these form modules are not isomorphic to each
other.
\end{proof}

\begin{remark}\label{rmk-1}
It follows that the isomorphism class of the form module $W_l(m)$
over $\bar{\tF}_q$ remains as one isomorphism class over $\tF_q$
when $l=\frac{m}{2}$ and decomposes into two isomorphism classes
$W_l^0(m)$ and $W_l^\delta(m)$ over $\tF_q$ when $l\neq\frac{m}{2}$.
The isomorphism class of the form module $D(m)$ over $\bar{\tF}_q$
remains as one isomorphism class over $\tF_q$.
\end{remark}

\section{Nilpotent orbits over $\tF_q$} In this section we study
the nilpotent orbits in the orthogonal Lie algebras over $\tF_q$ by
extending the method in \cite{Hes}. Let $V$ be a non-degenerate form
space over $\bar{\tF}_q$. An isomorphism class of form modules on
$V$ over $\bar{\tF}_q$ may decompose into several isomorphism
classes over
 $\tF_q$.

\begin{proposition}\label{prop-2}
Let $W$ be a form module
$(\lambda_1)_{\chi(\lambda_1)}^{m_1}(\lambda_2)_{\chi(\lambda_2)}^{m_2}
\cdots(\lambda_s)_{\chi(\lambda_s)}^{m_s}$ on the form space $V$.

$\mathrm{(i)}$ Assume $V$ is defective. The isomorphism class of $W$
over $\bar{\tF_q}$ decomposes into at most $2^{n_1}$ isomorphism
classes over $\tF_q$, where $n_1$ is the cardinality of $\{1\leq
i\leq s-1|\chi(\lambda_i)+\chi(\lambda_{i+1})\leq \lambda_i,\
\chi(\lambda_i)\neq\lambda_i/2\}$.

$\mathrm{(ii)}$ Assume $V$ is non-defective. The isomorphism class
of $W$ over $\bar{\tF_q}$ decomposes into  at most $2^{n_2}$
isomorphism classes over $\tF_q$, where $n_2$ is the cardinality of
$\{1\leq i\leq s|\chi(\lambda_i)+\chi(\lambda_{i+1})\leq \lambda_i,\
\chi(\lambda_i)\neq\lambda_i/2\}$ (here define
$\chi(\lambda_{s+1})=0$).
\end{proposition}
Note that we have two types of non-defective form spaces of
dimension $2n$ over $\tF_q$, $V^+$ with a quadratic form of Witt
index $n$ and $V^{-}$ with a quadratic form of Witt index $n-1$. We
define $O^{+}_{2n}(\tF_q)$ (resp. $O^{-}_{2n}(\tF_q)$) to be
$O(V^+)$ (resp. $O(V^-)$) and $\mathfrak{o}^{+}_{2n}(\tF_q)$ (resp.
$\mathfrak{o}^{-}_{2n}(\tF_q)$) to be $\Lo(V^+)$ (resp. $\Lo(V^-)$).
Let $SO^+_{2n}(\tF_q)=O^{+}_{2n}(\tF_q)\cap SO_{2n}(\bar{\tF}_q)$. A
form module on $V^+$ (resp. $V^-$) has orthogonal decomposition
$W_{l_1}^{\epsilon_1}(\lambda_1)\oplus\cdots\oplus
W_{l_{k}}^{\epsilon_{k}}(\lambda_{k})$
with $\#\{1\leq i\leq k|\epsilon_i=\delta\}$ being even (resp. odd).

\begin{corollary}\label{coro-n}
$\mathrm{(i)}$ The nilpotent $O_{2n+1}(\bar{\tF}_q)$-orbit
$(\lambda_1)_{\chi(\lambda_1)}^{m_1}\cdots(\lambda_s)_{\chi(\lambda_s)}^{m_s}
$ in $\mathfrak{o}_{2n+1}(\bar{\tF}_q)$ decomposes into at most
$2^{n_1}$ $O_{2n+1}(\tF_q)$-orbits in $\mathfrak{o}_{2n+1}(\tF_q)$.

$\mathrm{(ii)}$ If $\chi(\lambda_i)=\lambda_i/2\ ,i=1,\ldots,s$, the
nilpotent $O_{2n}(\bar{\tF}_q)$-orbit
$(\lambda_1)_{\chi(\lambda_1)}^{m_1}\cdots(\lambda_s)_{\chi(\lambda_s)}^{m_s}$
in $\mathfrak{o}_{2n}(\bar{\tF}_q)$ remains one
$O_{2n}^+({\tF}_q)$-orbit in $\mathfrak{o}^+_{2n}(\tF_q)$;
otherwise, it decomposes into at most $2^{n_2-1}$
$O^+_{2n}({\tF}_q)$-orbits in $\mathfrak{o}^+_{2n}(\tF_q)$ and at
most $2^{n_2-1}$ $O^-_{2n}({\tF}_q)$-orbits in
$\mathfrak{o}^-_{2n}(\tF_q)$.

Here $n_1,n_2$ are as in Proposition \ref{prop-2}.
\end{corollary}
\begin{remark}
In Corollary \ref{coro-n} $\mathrm{(ii)}$, if
$\chi(\lambda_i)=\lambda_i/2\ ,i=1,\ldots,s$, then $n$ is even; if
$\chi(\lambda_i)\neq\lambda_i/2$ for some $i$, then $n_2\geq 1$.
\end{remark}

Before we prove Proposition \ref{prop-2}, we need the following
lemma.

\begin{lemma}\label{lem-2}
$\mathrm{(i)}$ Assume $ k\geq m$ and $l\geq m$. We have
$W_l^0(k)\oplus D(m)\cong W_l^\delta(k)\oplus D(m)$ if and only if
$l+m>k.$

$\mathrm{(ii)}$ Assume $m>k$. We have
$ D(m)\oplus W_{k}^0(k)\cong D(m)\oplus W_{k}^\delta(k)$.

$\mathrm{(iii)}$ Assume $l_1\geq
l_2,\lambda_1-l_1\geq\lambda_2-l_2$. If $l_1+l_2>\lambda_1$, then
$W_{l_1}^0(\lambda_1)\oplus W_{l_2}^0(\lambda_2)\cong
W_{l_1}^\delta(\lambda_1)\oplus W_{l_2}^\delta(\lambda_2)$ and
$W_{l_1}^0(\lambda_1)\oplus W_{l_2}^\delta(\lambda_2)\cong
W_{l_1}^\delta(\lambda_1)\oplus W_{l_2}^0(\lambda_2)$.

$\mathrm{(iv)}$ Assume $l_1\geq l_2,\lambda_1-l_1\geq\lambda_2-l_2$.
If $l_1+l_2\leq\lambda_1$, then
$W_{l_1}^{\epsilon_1}(\lambda_1)\oplus
W_{l_2}^{\epsilon_2}(\lambda_2)\cong
W_{l_1}^{\epsilon_1'}(\lambda_1)\oplus
W_{l_2}^{\epsilon_2'}(\lambda_2)$ if and only if
$\epsilon_1=\epsilon_1',\epsilon_2=\epsilon_2'$, where $\epsilon_i,
\epsilon_i'=0$ or $\delta$, $i=1,2$.

\end{lemma}
\begin{proof}
We only prove $\mathrm{(i)}$. $\mathrm{(ii)}$-(iv) are proved
similarly. Take $v_1,w_1$ and $v_2,w_2$ such that $W_{l}^0(k)\oplus
D(m)=Av_1\oplus Aw_1\oplus Av_2\oplus Aw_2$ and
$\psi(v_1)=t^{2-2l},\psi(w_1)=0,\varphi(v_1,w_1)=t^{1-k},\psi(v_2)=t^{2-2m},\psi(w_2)=0,
\varphi(v_2,w_2)=t^{2-m},\varphi(v_1,v_2)=\varphi(v_1,w_2)=\varphi(w_1,v_2)=\varphi(w_1,w_2)=0.
$ Similarly, take $v_1',w_1'$ and $v_2',w_2'$ such that
$W_{l}^\delta(k)\oplus D(m)=Av_1'\oplus Aw_1'\oplus Av_2'\oplus
Aw_2'$ and $\psi(v_1')=t^{2-2l},\psi(w_1')=\delta
t^{2l-2k},\varphi(v_1',w_1')=t^{1-k},\psi(v_2')=t^{2-2m},\psi(w_2')=0,
\varphi(v_2',w_2')=t^{2-m},\varphi(v_1',v_2')=\varphi(v_1',w_2')=\varphi(w_1',v_2')=\varphi(w_1',w_2')=0.
$

We have $W_{l}^0(k)\oplus D(m)\cong W_{l}^\delta(k)\oplus D(m)$ if
and only if there exists an $A$-module isomorphism $g:V\rightarrow
V$ such that $\psi(gv)=\psi(v)$ and $\varphi(gv,gw)=\varphi(v,w)$
for any $v,w\in V$. Assume
\begin{eqnarray*}
&&gv_j=\sum\limits_{i=0}^{k-1}(a_{j,i}t^iv_1'+b_{j,i}t^iw_1')+\sum\limits_{i=0}^{m-1}c_{j,i}t^{i}v_{2}'
+\sum\limits_{i=0}^{m-2}d_{j,i}t^{i}w_{2}',\\
&&gw_j=\sum\limits_{i=0}^{k-1}(e_{j,i}t^iv_1'+f_{j,i}t^iw_1')+\sum\limits_{i=0}^{m-1}g_{j,i}t^{i}v_{2}'
+\sum\limits_{i=0}^{m-2}h_{j,i}t^{i}w_{2}',j=1,2.
\end{eqnarray*}
Then $W_l^0(k)\oplus D(m)\cong W_l^\delta(k)\oplus D(m)$ if and only
if the equations $\psi(gv_i)=\psi(v_i)$,
$\psi(gw_i)=\psi(w_i),\varphi(gv_i,gv_j)=\varphi(v_i,v_j),
\varphi(gv_i,gw_j)=\varphi(v_i,w_j),\varphi(gw_i,gw_j)=\varphi(w_i,w_j)$
have solutions.

If $l+m\leq k$, some equations are
$e_{1,2l-k-1}^2+e_{1,2l-k-1}=\delta$ $(l\neq\frac{k+1}{2})$ or
$a_{1,0}^2+a_{1,0}b_{1,0}=1,
e_{1,0}^2+e_{1,0}f_{1,0}=\delta,a_{1,0}f_{1,0}+b_{1,0}e_{1,0}=1$
$(l=\frac{k+1}{2})$. By the definition of $\delta$, this has no
solution for $e_{1,2l-k-1}$ or $a_{1,0},b_{1,0},e_{1,0},f_{1,0}$,
which implies that $W_l^0(k)\oplus D(m)\ncong W_l^\delta(k)\oplus
D(m)$.

If $l+m>k$, let
$gv_1=v_1',gw_1=w_1'+\sqrt{\delta}t^{l+m-k-1}v_2',gv_2=v_2',gw_2=w_2'+\sqrt{\delta}t^{l}v_1'$.
This is a solution for the equations. It follows that
$W_l^0(k)\oplus D(m)\cong W_l^\delta(k)\oplus D(m)$.
\end{proof}

\begin{proof}[Proof of Proposition \ref{prop-2}]
We prove (i). One can prove (ii) similarly. By the classification in
subsection \ref{sec-1}, we can rewrite the symbol as
$$(\lambda_1)_{\chi(\lambda_1)}^{2}(\lambda_2)_{\chi(\lambda_2)}^{2}\cdots(\lambda_k)_{\chi(\lambda_{k})}^{2}
(\lambda_{k+1})_{\lambda_{k+1}}(\lambda_{k+1}-1)_{(\lambda_{k+1}-1)}(\lambda_{k+2})_{\lambda_{k+2}}^{2}\cdots
(\lambda_{k+l})_{\lambda_{k+l}}^{2},$$ where
$\lambda_1\geq\cdots\geq\lambda_{k+1}>\lambda_{k+2}\geq\cdots\geq\lambda_{k+l}$,
$\chi(\lambda_i)\geq\chi(\lambda_{i+1})$ and
$\lambda_i-\chi(\lambda_i)\geq \lambda_{i+1}-\chi(\lambda_{i+1})$
for $1\leq i\leq k$ (by abuse of notation, we still use $\lambda$),
$l\geq 1$. A representative $W$ for this isomorphism class over
$\bar{\tF}_q$ is $ W_{\chi(\lambda_1)}(\lambda_1)\oplus\cdots\oplus
W_{\chi(\lambda_k)}(\lambda_{k})\oplus D(\lambda_{k+1})\oplus
W_{\lambda_{k+2}}(\lambda_{k+2})\oplus\cdots\oplus
W_{\lambda_{k+l}}(\lambda_{k+l})$. By Proposition \ref{prop-nind}
and Remark \ref{rmk-1}, in order to study the isomorphism classes
into which the isomorphism class of $W$ over $\bar{\tF}_q$
decomposes over $\tF_q$, it is enough to study the isomorphism
classes of form modules of the form $
W_{\chi(\lambda_1)}^{\epsilon_1}(\lambda_1)\oplus\cdots\oplus
W_{\chi(\lambda_k)}^{\epsilon_k}(\lambda_{k})\oplus
D(\lambda_{k+1})\oplus
W_{\lambda_{k+2}}^{\epsilon_{k+2}}(\lambda_{k+2})\oplus\cdots\oplus
W_{\lambda_{k+l}}^{\epsilon_{k+l}}(\lambda_{k+l})$, where
$\epsilon_i=0$ or $\delta$. Thus it suffices to show that modules of
the above form have at most $2^{n_1}$ isomorphism classes.

We have $n_1=\#\{1\leq i\leq
k|\chi(\lambda_i)+\chi(\lambda_{i+1})\leq \lambda_i,\
\chi(\lambda_i)\neq\lambda_i/2\}$. Suppose $i_1,i_2,\ldots,i_{n_1}$
are such that $1\leq i_j\leq
k,\chi(\lambda_{i_j})+\chi(\lambda_{i_j+1})\leq \lambda_{i_j},\
\chi(\lambda_{i_j})\neq\lambda_{i_j}/2,j=1,\ldots,n_1$. Then using
Lemma \ref{lem-2} one can easily show that a module of the above
form is isomorphic to one of the following modules:
$V_1^{\epsilon_1}\oplus\cdots\oplus V_{n_1}^{\epsilon_{n_1}}\oplus
V_{n_1+1}$, where
$V_{t}^{\epsilon_t}=W_{\chi(\lambda_{i_{t-1}+1})}^0(\lambda_{i_{t-1}+1})\oplus\cdots\oplus
W_{\chi(\lambda_{i_{t}-1})}^0(\lambda_{i_{t}-1})\oplus
W_{\chi(\lambda_{i_{t}})}^{\epsilon_t}(\lambda_{i_{t}})$,
$t=1,\ldots,n_1$, $i_0=0$, $\epsilon_t=0$ or $\delta$ and
$V_{n_1+1}=W_{\chi(\lambda_{i_{n_1}+1})}^0(\lambda_{i_{n_1}+1})\oplus\cdots\oplus
W_{\chi(\lambda_k)}^0(\lambda_k)\oplus D(\lambda_{k+1})\oplus
W_{\lambda_{k+2}}^0(\lambda_{k+2})\oplus\cdots\oplus
W_{\lambda_{k+l}}^0(\lambda_{k+l})$. Thus $\mathrm{(i)}$ is proved.
\end{proof}

\section{Number of nilpotent orbits over $\tF_q$}

 \subsection{} In this subsection we
recall Spaltenstein's parametrization of nilpotent orbits by pairs
of partitions in $\Lo(\bar{\tF}_q)$ (see \cite{Spal}).

For $\mathfrak{o}_{2n+1}(\bar{\tF}_q)$, the orbit
$$(\lambda_1)_{\chi(\lambda_1)}^{2}\cdots(\lambda_k)_{\chi(\lambda_{k})}^{2}
(\lambda_{k+1})_{\lambda_{k+1}}(\lambda_{k+1}-1)_{(\lambda_{k+1}-1)}
(\lambda_{k+2})_{\lambda_{k+2}}^{2}\cdots
(\lambda_{k+l})_{\lambda_{k+l}}^{2}$$ is written as
$$(\alpha_1+\beta_1)_{(\alpha_1+1)}^{2}\cdots(\alpha_k+\beta_k)_{(\alpha_k+1)}^{2}
(\alpha_{k+1}+1)_{(\alpha_{k+1}+1)}(\alpha_{k+1})_{\alpha_{k+1}}\\
(\alpha_{k+2})_{\alpha_{k+2}}^{2}\cdots
(\alpha_{k+l})_{\alpha_{k+l}}^{2}$$  and the corresponding pair of
partitions is $(\alpha,\beta)$, where
$\alpha=(\alpha_1\ldots,\alpha_{k+l})$ and
$\beta=(\beta_1,\ldots,\beta_k)$ satisfy
$\alpha_1\geq\cdots\geq\alpha_{k+l}\geq 0,\
\beta_1\geq\cdots\geq\beta_k\geq 1$ and $|\alpha|+|\beta|=n$.

For $\mathfrak{o}_{2n}(\bar{\tF}_q)$, the orbit
$$(\lambda_1)_{\chi(\lambda_1)}^{2}\cdots(\lambda_k)_{\chi(\lambda_{k})}^{2}
$$ is written as $$(\alpha_1+\beta_1)_{\alpha_1}^{2}\cdots
(\alpha_k+\beta_k)_{\alpha_k}^{2}$$ and the corresponding pair of
partitions is $(\alpha,\beta)$, where
$\alpha=(\alpha_1,\ldots,\alpha_k)$ and $
\beta=(\beta_1,\ldots,\beta_k)$ satisfy
$\alpha_1\geq\cdots\geq\alpha_k\geq 1,\
\beta_1\geq\cdots\geq\beta_k\geq 0$ and $|\alpha|+|\beta|=n$.

\subsection{} In this subsection we study the number of nilpotent
orbits over $\tF_q$. Denote by $p_2(n)$ the cardinality of the set
of pairs of partitions $(\alpha,\beta)$ such that
$|\alpha|+|\beta|=n$ and $p(k)$ the number of partitions of the
integer $k$.
\begin{proposition}\label{prop-num}
$\mathrm{(i)}$ The number of nilpotent $O_{2n+1}(\tF_q)$-orbits in
$\mathfrak{o}_{2n+1}(\tF_q)$ is at most $p_2(n)$.

$\mathrm{(ii)}$ The number of nilpotent $O_{2n}^+(\tF_q)$-orbits in
$\mathfrak{o}_{2n}^+(\tF_q)$
 is at most $\frac{1}{2}p_2(n)$ if $n$ is odd and is
at most $\frac{1}{2}(p_2(n)+p(\frac{n}{2}))$ if $n$ is even.

\end{proposition}
\begin{proof}
$\mathrm{(i)}$ The set of nilpotent orbits in
$\mathfrak{o}_{2n+1}(\bar{\tF}_q)$ is mapped bijectively to the set
$\{(\alpha,\beta)||\alpha|+|\beta|=n,\beta_i\leq
\alpha_i+2\}:=\Delta$(\cite{Spal}). By Corollary \ref{coro-n} (i), a
nilpotent orbit in $\mathfrak{o}_{2n+1}(\bar{\tF}_q)$ corresponding
to $(\alpha,\beta)\in\Delta,\alpha=(\alpha_1,\ldots,\alpha_s),
\beta=(\beta_1,\ldots,\beta_t)$ splits into at most $2^{n_1}$ orbits
in $\mathfrak{o}_{2n+1}(\tF_q)$, where $n_1=\#\{1\leq i\leq
t|\alpha_{i+1}+2\leq\beta_i<\alpha_i+2\}$. We associate to this
orbit $2^{n_1}$ pairs of partitions as follows. Suppose
$r_1,r_2,\ldots,r_{n_1}$ are such that
$\alpha_{r_{i}+1}+2\leq\beta_{r_i}<\alpha_{r_i}+2,i=1,\ldots,n_1$.
Let
\begin{eqnarray*}&&\alpha^{1,i}=(\alpha_{r_{i-1}+1},\ldots,\alpha_{r_i}),
\beta^{1,i}=(\beta_{r_{i-1}+1},\ldots,\beta_{r_i}),\\
&&\alpha^{2,i}=(\beta_{r_{i-1}+1}-2,\ldots,\beta_{r_i}-2),
\beta^{2,i}=(\alpha_{r_{i-1}+1}+2,\ldots,\alpha_{r_i}+2),i=1,\ldots,n_1,\\
&&\alpha^{n_1+1}=(\alpha_{r_{n_1}+1},\ldots,\alpha_{s}),
\beta^{n_1+1}=(\beta_{r_{n_1}+1},\ldots,\beta_{t}). \end{eqnarray*}
We associate to $(\alpha,\beta)$ the pairs of partitions
$(\tilde{\alpha}^{\epsilon_1,\ldots,\epsilon_{n_1}},
\tilde{\beta}^{\epsilon_1,\ldots,\epsilon_{n_1}})$,
$$\tilde{\alpha}^{\epsilon_1,\ldots,\epsilon_{n_1}}=
(\alpha^{\epsilon_1,1},\ldots,\alpha^{\epsilon_{n_1},n_1},\alpha^{n_1+1}),
\tilde{\beta}^{\epsilon_1,\ldots,\epsilon_{n_1}}=
(\beta^{\epsilon_1,1},\ldots,\beta^{\epsilon_{n_1},{n_1}},\beta^{{n_1}+1}),$$
where $\epsilon_i\in\{1,2\},i=1,\ldots,{n_1}$.

Notice that the pairs of partitions
$(\tilde{\alpha}^{\epsilon_1,\ldots,\epsilon_{n_1}},
\tilde{\beta}^{\epsilon_1,\ldots,\epsilon_{n_1}})$ are distinct and
among them only $(\alpha,\beta)=(\tilde{\alpha}^{1,\ldots,1},
\tilde{\beta}^{1,\ldots,1})$ is in $\Delta$. One can verify that the
set of all pairs of partitions
$(\tilde{\alpha}^{\epsilon_1,\ldots,\epsilon_{n_1}},
\tilde{\beta}^{\epsilon_1,\ldots,\epsilon_{n_1}})$ constructed as
above for $(\alpha,\beta)\in\Delta$ is equal to the set
$\{(\alpha,\beta)||\alpha|+|\beta|=n\}$, which has cardinality
$p_2(n)$. But the number of nilpotent orbits in
$\mathfrak{o}_{2n+1}(\tF_q)$ is no greater than the cardinality of
the former set. Thus $\mathrm{(i)}$ is proved.

$\mathrm{(ii)}$ Similarly, the set of nilpotent orbits in
$\mathfrak{o}_{2n}(\bar{\tF}_q)$ is mapped bijectively to the set
$\{(\alpha,\beta)||\alpha|+|\beta|=n,\beta_i\leq
\alpha_i\}:=\Delta'$(\cite{Spal}). By Corollary \ref{coro-n} (ii), a
nilpotent orbit in $\mathfrak{o}_{2n}(\bar{\tF}_q)$ corresponding to
$(\alpha,\beta)\in\Delta'$ with
$\alpha=(\alpha_1,\alpha_2,\ldots,\alpha_s),\beta=(\beta_1,\beta_2,\ldots,\beta_s)$
and $\alpha\neq\beta$ splits into at most $2^{n_2-1}$ orbits in
$\mathfrak{o}^{+}_{2n}(\tF_q)$, where
$n_2=\#\{i|\alpha_{i+1}\leq\beta_i<\alpha_i\}-1$. We associate to
this orbit $2^{n_2-1}$ pairs of partitions as follows. We can assume
$\alpha_{s}\neq 0$. Suppose $r_{1},r_{2},\ldots,r_{{n_2}}$ are such
that
$\alpha_{r_{i}+1}\leq\beta_{r_{i}}<\alpha_{r_{i}},i=1,\ldots,{n_2}$.
Let
\begin{eqnarray*}
&&\alpha^{1,i}=(\alpha_{r_{i-1}+1},\ldots,\alpha_{r_i}),
\beta^{1,i}=(\beta_{r_{i-1}+1},\ldots,\beta_{r_i}),\\
&&\alpha^{2,i}=(\beta_{r_{i-1}+1},\ldots,\beta_{r_i}),
\beta^{2,i}=(\alpha_{r_{i-1}+1},\ldots,\alpha_{r_i}),i=1,\ldots,{n_2},\\
&&\alpha^{{n_2}+1}=(\alpha_{r_{{n_2}}+1},\ldots,\alpha_{s}),
\beta^{{n_2}+1}=(\beta_{r_{{n_2}}+1},\ldots,\beta_{s}).
\end{eqnarray*}
We have $2^{n_2}$ distinct pairs of partitions
$(\tilde{\alpha}^{\epsilon_1,\ldots,\epsilon_{n_2}},\tilde{\beta}^{\epsilon_1,\ldots,\epsilon_{n_2}}),$$$
\tilde{\alpha}^{\epsilon_1,\ldots,\epsilon_{n_2}}=(\alpha^{\epsilon_1,1},\ldots,\\
\alpha^{\epsilon_{n_2},{n_2}},\alpha^{{n_2}+1}),
\tilde{\beta}^{\epsilon_1,\ldots,\epsilon_{n_2}}=(\beta^{\epsilon_1,1},\ldots,\beta^{\epsilon_{n_2},{n_2}},\beta^{{n_2}+1}),$$
where $\epsilon_i\in\{1,2\},i=1,\ldots,{n_2}$. We show that in these
pairs of partitions $(\alpha',\beta')$ appears if and only if
$(\beta',\alpha')$ appears. In fact we have $\alpha_{i}=\beta_{i}$,
for $i>r_{{n_2}}$, which implies $\alpha^{{n_2}+1}=\beta^{{n_2}+1}$.
Thus we have $$(\tilde{\alpha}^{\epsilon_1+1(mod
2),\ldots,\epsilon_{n_2}+1(mod 2)}, \tilde{\beta}^{\epsilon_1+1(mod
2),\ldots,\epsilon_{n_2}+1(mod 2)})
=(\tilde{\beta}^{\epsilon_1,\ldots,\epsilon_{n_2}},\tilde{\alpha}^{\epsilon_1,\ldots,\epsilon_{n_2}}).$$
Hence we can identify $(\alpha',\beta')$ with $(\beta',\alpha')$,
and then associate $2^{{n_2}-1}$ pairs of partitions to the
nilpotent orbit corresponding to $(\alpha,\beta)$.

One can verify that the set of all pairs of partitions associated to
 $(\alpha,\beta)\in\Delta'$ as above is in bijection with the set of
pairs of partitions $(\alpha,\beta)$ such that $|\alpha|+|\beta|=n$
with $(\alpha,\beta)$ identified with $(\beta,\alpha)$, which has
cardinality $\frac{1}{2}p_2(n)$ if $n$ is odd and
$\frac{1}{2}(p_2(n)+p(\frac{n}{2}))$ if $n$ is even. Thus (ii)
follows.
\end{proof}

\begin{corollary}\label{coro}
The number of nilpotent $SO_{2n}^+(\tF_q)$-orbits in
$\mathfrak{o}_{2n}^+(\tF_q)$ is at most $\frac{1}{2}p_2(n)$ if $n$
is odd and is at most $\frac{1}{2}(p_2(n)+3p(\frac{n}{2}))$ if $n$
is even.
\end{corollary}
\begin{proof}
We show that the $O_{2n}^+(\tF_q)$ orbits that split into two
$SO_{2n}^+(\tF_q)$-orbits are exactly the orbits corresponding to
the pairs of partitions of the form $(\alpha,\alpha)$. The number of
these orbits is $p(\frac{n}{2})$. Let $x$ be a nilpotent element in
$\mathfrak{o}_{2n}^+(\tF_q)$. The $O^+_{2n}(\tF_q)$-orbit of $x$
splits into two $SO^+_{2n}(\tF_q)$-orbits if and only if the
centralizer $Z_{O^+_{2n}(\tF_q)}(x)\subset SO^+_{2n}(\tF_q)$. It is
enough to show that for an indecomposable module $V$,
$Z_{O^+_{2n}(\tF_q)}(V)\subset SO^+_{2n}(\tF_q)$ if and only if
$\chi(m)\leq\frac{1}{2}m$, for all $m\in\mathbb{N}$.

Assume $V=W_l^0(m)$ or $W_l^\delta(m)$, $l\geq\frac{m+1}{2}$. Let
$\epsilon=0$ or $\delta$ and $v_1^\epsilon,v_2^\epsilon$ be such
that $W_l^{\epsilon}(m)=A v_1^\epsilon\oplus A v_2^\epsilon$ and
$\psi(v_1^\epsilon)=t^{2-2l}$, $\psi(v_2^\epsilon)=\epsilon
t^{-2m+2l}$, $\varphi(v_1^\epsilon,v_2^\epsilon)=t^{1-m}$. Let
$w_2^\epsilon=v_2^\epsilon+t^{2l-1-m}v_1^\epsilon$. Define
 $u^\epsilon$ by $u^\epsilon(a_1v_1^\epsilon+a_2v_2^\epsilon)=a_1v_1^\epsilon+a_2w_2^\epsilon$.
 Then $u^\epsilon\in Z_{O^+_{2n}(\tF_q)}(W_l^{\epsilon}(m))$,
 but $u^\epsilon\notin SO^+_{2n}(\tF_q)$ (see \cite{Hes}). This
 shows that $Z_{O^+_{2n}(\tF_q)}(V)\nsubseteq SO^+_{2n}(\tF_q)$.

Assume $V=W_l^0(m)$, $l=m/2$. Let $v_1,v_2$ be such that
$W_l^{0}(m)=Av_1\oplus Av_2$ and $\psi(v_1)=t^{2-m}$, $\psi(v_2)=0$,
$\varphi(v_1,v_2)=t^{1-m}$. Let $W$ be the subspace of $V$ spanned
by $t^{\frac{m}{2}}v_1,t^{\frac{m}{2}+1}v_1,
\ldots,t^{m-1}v_1,t^{\frac{m}{2}}v_2$, $t^{\frac{m}{2}+1}v_2,
\ldots,t^{m-1}v_2$. Then $W$ is a maximal totally singular subspace
in $V$ and $Z_{O^+_{2n}(\tF_q)}(V)$ stabilizes $W$. Hence
$Z_{O^+_{2n}(\tF_q)}(V)\subset SO^+_{2n}(\tF_q).$
\end{proof}

\section{Springer corresopondence}

Throughout subsections \ref{ssec-2}-\ref{ssec-3}, let $G$ be a
connected adjoint algebraic group of type $B_r,C_r$ or $D_r$ over
$\tk$ and $\mathfrak{g}$ be the Lie algebra of $G$. Fix a Borel
subgroup $B$ of $G$ with Levi decomposition $B=TU$. Let
$\mathfrak{b},\mathfrak{t}$ and $\mathfrak{n}$ be the Lie algebra of
$B,T$ and $U$ respectively. Let $\cB$ be the variety of Borel
subgroups of $G$.

We construct in this section the Springer correspondence for the Lie
algebra $\Lg$ following \cite{Lu1,Lu4}. The construction and proofs
are essentially the same as (actually simpler than) those for the
unipotent case in \cite{Lu1,Lu4}. For completeness, we include the
proofs here. We only point out that Lemma 6.2 is essential for the
construction. The lemma is probably well-known for which we include
an elementary proof.

\subsection{}\label{ssec-2}
Let $Z$ be the variety $\{(x,B_1,B_2)\in
\mathfrak{g}\times\cB\times\cB|x\in\Lb_1\cap\Lb_2\}$ and $ Z'$ be
the Steinberg variety \cite{st2} $\{(x,B_1,B_2)\in
\mathfrak{g}\times\cB\times\cB|x\in\Ln_1\cap\Ln_2\}. $ Denote $r$
the dimension of $T$. Let $\rc$ be a nilpotent orbit in $\Lg$. A
stronger version of the following lemma in the group case is due to
Springer, Steinberg and Spaltenstein (see for example \cite{Spa3}).
We include a proof for the Lie algebra case here.
\begin{lemma}\label{prop-dim}
$\mathrm{(i)}$ We have $\dim(\rc\cap\Ln)\leq\frac{1}{2}\dim \rc$.

$\mathrm{(ii)}$ Given $x\in \rc$, we have
$\dim\{B_1\in\cB|x\in\Ln_1\}\leq(\dim G-r-\dim \rc)/2$.

$\mathrm{(iii)}$ We have $\dim Z=\dim G$ and $\dim Z'=\dim G-r$.
\end{lemma}
\begin{proof}
We have a partition
$Z=\cup_{\mathcal{\mathcal{O}}}Z_{\mathcal{\mathcal{O}}}$ according
to the $G$-orbits ${\mathcal{\mathcal{O}}}$ on $\cB\times\cB$ where
$Z_\mathcal{O}=\{(x,B_1,B_2)\in Z|(B_1,B_2)\in\mathcal{O}\}$. Define
in the same way a partition
$Z'=\cup_{\mathcal{\mathcal{O}}}Z'_{\mathcal{\mathcal{O}}}$.
Consider the maps from $Z_{\cO}$ and $Z'_{\cO}$ to $\cO$:
$(x,B_1,B_2)\mapsto (B_1,B_2)$. We have $\dim Z_{\cO}=\dim
(\Lb_1\cap\Lb_2)+\dim\cO=\dim G$ and $\dim Z_{\cO}'=\dim
(\Ln_1\cap\Ln_2)+\dim\cO=\dim G-r$. Thus $\mathrm{(iii)}$ follows.

Let $Z'(\rc)=\{(x,B_1,B_2)\in Z'|x\in \rc\}\subset Z'$. From (iii),
we have $\dim Z'(\rc)\leq \dim G-r$. Consider the map
$Z'(\rc)\rightarrow \rc,\ (x,B_1,B_2)\mapsto x$. We have $\dim
Z'(\rc)=\dim\rc+2\dim\{B_1\in\cB|x\in\Ln_1\}\leq \dim G-r$. Thus
(ii) follows.

Consider the variety $\{(x,B_1)\in \rc\times\cB|x\in\Ln_1\}$. By
projecting it to the  first coordinate, and using (ii), we see that
it has dimension $\leq(\dim G-r+\dim \rc)/2$. If we project it to
the second coordinate, we get $\dim(\rc\cap\Ln)+\dim\cB\leq(\dim G-r
+\dim \rc)/2$ and (i) follows.
\end{proof}

\subsection{}
Recall that an element $x$ in $\Lg$ is called regular if the
connected centralizer $Z_G^0(x)$ in $G$ is a maximal torus of $G$.
\begin{lemma}\label{ssreg}
There exist regular semisimple elements in $\Lg$ and they form an
open dense subset in $\Lg$.
\end{lemma}
\begin{proof}
We first show that regular semisimple elements exists in $\Lg$.

 (i)
$G=SO(2n+1)$. Let $V$ be a $2n+1$ dimensional vector space equipped
with a non-degenerate quadratic form $Q:V\rightarrow k$. Let
$\langle\ ,\ \rangle$ be the bilinear form associated to $Q$ defined
by $\langle v,w\rangle=Q(v+w)+Q(v)+Q(w)$. Then $G$ is defined as
$\{g\in GL(V)|\ Q(gv)=Q(v),\ \forall\ v\in V\}$. Since $Q$ is
non-degenerate, we have $\dim V^{\perp}=1$. Let
$V^{\perp}=\text{span}\{v_0\}$, we have $gv_0=v_0$, for any $g\in
G$. Hence the Lie algebra is $\{x\in \End(V)|\langle
xv,v\rangle=0,\forall\ v\in V; xv_0=0\}$. With respect to a suitable
basis, we can assume $Q(v)=v^tBv$, where $B=\left(
\begin{array}{ccc}
0 & I_{n\times n} & 0 \\   0_{n\times n} & 0 & 0 \\
0 & 0 & 1 \\
\end{array}
\right).$ A maximal torus of $G$ is
$T=\{\diag(t_1,t_2,\ldots,t_n,1/t_1,1/t_2,\ldots 1/t_n,1)|t_i\in
k^*,i=1,\ldots,n\}.$ The Lie algebra of $T$ is
$\Lt=\{\diag(x_1,x_2,\ldots,x_n,x_1,x_2,\ldots x_n,0)|x_i\in
k,i=1,\ldots,n\}.$ Since every semisimple element in $\Lg$ is
conjugate to an element in $\Lt$ under the adjoint action of $G$, it
is enough to consider elements in $\Lt$.

Let $x=\diag(x_1,x_2,\ldots,x_n,x_1,x_2,\ldots x_n,0)$, where
$x_i\neq x_j$, for any $i\neq j$ and $x_i\neq 0$ for any $i$ (such
$x$ exists). It can be easily verified that $Z_G(x)$ consists of
elements of the form $g=\left(
\begin{array}{ccc}
A_1 & A_2 & 0 \\   A_3 & A_4 & 0 \\
0 & 0 & 1 \\
\end{array}
\right)$ where $A_i=\diag(a_i^1,a_i^2,\ldots,a_i^n),i=1,2,3,4,$
 satisfy $a_1^ja_3^j=a_2^ja_4^j=0$ and
$a_1^ja_4^j+a_2^ja_3^j=1,j=1,\ldots,n.$ Hence we see that
$Z^0_G(x)=T$ and $x$ is regular.

(ii) $G$ is the adjoint group of type $C_n$. We have the following
construction of $G$. Let $V$ be a $2n$ dimensional vector space
equipped with a non-degenerate symplectic form $\langle\ ,\
\rangle:V\times V\rightarrow k$. Then $G$ is defined as
$$\{(g,\lambda)\in GL(V)\times k^*|\forall\ v,w\in V, \langle
gv,gw\rangle=\lambda \langle v,w\rangle\}/\{(\mu I,\mu^2)|\mu\in
k^*\}.$$ Hence the Lie algebra $\Lg$ is $$\{(x,\lambda)\in
\End(V)\times k|\forall\ v,w\in V, \langle xv,w\rangle+\langle
v,xw\rangle=\lambda \langle v,w\rangle\}/\{(\mu I,0)|\mu\in k\}.$$
With respect to a suitable basis, we can assume $\langle
v,w\rangle=v^tAw$ where $A=\left(
\begin{array}{cc}
0 & I_{n\times n}  \\   I_{n\times n} & 0
\end{array}
\right).$ A maximal torus of $G$ is
$T=\{\diag(t_1,t_2,\ldots,t_n,\lambda/t_1,\lambda/t_2,\ldots
\lambda/t_n)|t_i\in k^*,i=1,\ldots,n;\lambda\in k^*\}/\{(\mu
I,\mu^2)|\mu\in k^*\}.$ The Lie algebra of $T$ is
$\Lt=\{\diag(x_1,x_2,\ldots,x_n,\lambda+x_1,\lambda+x_2,\ldots
\lambda+x_n)|x_i\in k,i=1,\ldots,n;\lambda\in k\}/\{(\mu I,0)|\mu\in
k\}.$

Let
$x=\overline{\diag(x_1,x_2,\ldots,x_n,x_1+\lambda,x_2+\lambda,\ldots
x_n+\lambda)}\in\Lt$, where $x_i\neq x_j$, for any $i\neq j$ and
$\lambda\neq x_i+x_j$ for any $i,j$(such $x$ exists). It can be
verified that the elements in $Z_G(x)$ must have the form $g=\left(
\begin{array}{cc}
A_1 & A_2  \\   A_3 & A_4
\end{array}
\right)$ where each row and each column of $G$ has only one nonzero
entry and $(A_1)_{ij}\neq 0\Leftrightarrow(A_4)_{ij}\neq
0,(A_2)_{ij}\neq 0\Leftrightarrow(A_3)_{ij}\neq 0$. Hence we see
that $Z^0_G(x)=T$ and $x$ is regular.

(iii) $G$ is the adjoint group of type $D_n$. We have the following
construction of $G$. Let $V$ be a $2n$ dimensional vector space
equipped with a non-defective quadratic form $Q:V\rightarrow k$. Let
$\langle\ ,\ \rangle$ be the bilinear form associated to $Q$ defined
by $\langle v,w\rangle=Q(v+w)+Q(v)+Q(w)$. Then $G$ is defined as
$$\{(g,\lambda)\in GL(V)\times k^*|\forall\ v\in V, Q(gv)=\lambda
Q(v)\}/\{(\mu I,\mu^2)|\mu\in k^*\}.$$ Hence the Lie algebra $\Lg$
is $$\{(x,\lambda)\in \End(V)\times k|\forall\ v\in V, \langle
xv,v\rangle=\lambda Q(v)\}/\{(\mu I,0)|\mu\in k\}.$$ With respect to
a suitable basis, we can assume $Q(v)=v^tBv$ where $B=\left(
\begin{array}{cc}
0 & I_{n\times n}  \\   0 & 0
\end{array}
\right).$ A maximal torus of $G$ is
$T=\{\diag(t_1,t_2,\ldots,t_n,\lambda/t_1,\lambda/t_2,\ldots
\lambda/t_n)|t_i\in k^*,i=1,\ldots,n;\lambda\in k^*\}/\{(\mu
I,\mu^2)|\mu\in k^*\}.$ The Lie algebra of $T$ is
$\Lt=\{\diag(x_1,x_2,\ldots,x_n,\lambda+x_1,\lambda+x_2,\ldots
\lambda+x_n)|x_i\in k,i=1,\ldots,n;\lambda\in k\}/\{(\mu I,0)|\mu\in
k\}.$

Let
$x=\overline{\diag(x_1,x_2,\ldots,x_n,x_1+\lambda,x_2+\lambda,\ldots
x_n+\lambda)}\in\Lt$ where $x_i\neq x_j$, for any $i\neq j$ and
$\lambda\neq x_i+x_j$ for any $i\neq j$ (such $x$ exists), then
similarly one can show that $Z^0_G(x)=T$ and $x$ is regular.

Now denote by $\Lt_0$ the set of regular elements in $\Lt$. From the
above construction, one easily see that $\Lt_0$ is a dense subset in
$\Lt$. Thus $\dim\Lt_0=\dim\Lt=r$. Consider the map $\Lt_0\times
G\rightarrow \Lg,\ (x,g)\mapsto \Ad(g)x$. The fiber at $y$ in the
image is $\{(x,g)\in\Lt_0\times G|\Ad(g)x=y\}$. We consider the
projection of $\{(x,g)\in\Lt_0\times G|\Ad(g)x=y\}$ to the first
coordinate. The fiber of this projection at $x\in\Lt_0$ is
isomorphic to $Z_G(x)$, which has dimension $r$, and the image of
this projection is finite. Hence $\dim\{(x,g)\in\Lt_0\times
G|\Ad(g)x=y\}=r$. It follows that the image of the map $\Lt_0\times
G\rightarrow \Lg,\ (x,g)\mapsto \Ad(g)x$ has dimension equal to
$\dim(\Lt_0\times G)-r=\dim\Lg$. This proves the lemma.
\end{proof}

\begin{remark}
Lemma \ref{ssreg} is not always true when $G$ is not adjoint.
\end{remark}

\subsection{}
Let $Y$ be the set of semisimple regular elements in $\Lg$. By Lemma
\ref{ssreg}, $\dim Y=\dim G$. Let $\widetilde{Y}=\{(x,gT)\in Y\times
G/T|\Ad(g^{-1})(x)\in \Lt_0\}$. Define $\pi:
\widetilde{Y}\rightarrow Y$ by $\pi(x,gT)=x$. The Weyl group
$W=NT/T$ acts (freely) on $\widetilde{Y}$ by
$n:(x,gT)\mapsto(x,gn^{-1}T)$.
\begin{lemma}\label{lem-3}
$\pi: \widetilde{Y}\rightarrow Y$ is a principal $W$-bundle.
\end{lemma}
\begin{proof}
We show that if $x\in \Lg, g,g'\in G$ are such that
$\Ad(g^{-1})x\in\Lt_0$ and $\Ad(g'^{-1})x\in\Lt_0$, then
$g'=gn^{-1}$ for some $n\in NT$. Let
$\Ad(g^{-1})x=t_1\in\Lt_0,\Ad(g'^{-1})x=t_2\in\Lt_0$, then we have
$Z_G^0(x)=Z_G^0(\Ad(g)t_1)=gZ_G^0(t_1)g^{-1}=gTg^{-1}$, similarly
$Z_G^0(x)=Z_G^0(\Ad(g')t_2)=g'Z_G^0(t_2)g'^{-1}=g'Tg'^{-1}$, hence
$g'^{-1}g\in NT$.
\end{proof}

Let $X=\{(x,gB)\in\Lg\times G/B|\Ad(g^{-1})x\in\Lb\}$. Define
$\varphi:X\rightarrow\Lg$ by $\varphi(x,gB)=x$. The map $\varphi$ is
$G$-equivariant with $G$-action on $X$ given by
$g_0:(x,gB)\mapsto(\Ad(g_0)x,g_0gB)$.
\begin{lemma}
$\mathrm{(i)}$ $X$ is an irreducible variety of dimension equal to
$\dim G$.

$\mathrm{(ii)}$ $\varphi$ is proper and $\varphi(X)=\Lg=\bar{Y}$.

$\mathrm{(iii)}$ $(x,gT)\rightarrow(x,gB)$ is an isomorphism
$\gamma:\widetilde{Y}\xrightarrow{\sim}\varphi^{-1}(Y)$.
\end{lemma}
\begin{proof}
(i) and (ii) are easy. For (iii), we only prove that $\gamma$ is a
bijection. First we show that $\gamma$ is injective. Suppose
$(x_1,g_1T),(x_2,g_2T)\in\widetilde{Y}$ are such that
$(x_1,g_1B)=(x_2,g_2B)$, then we have
$\Ad(g_1^{-1})(x_1)\in\Lt_0,\Ad(g_2^{-1})(x_2)\in\Lt_0$ and
$x_1=x_2,g_2^{-1}g_1\in B$. Similar argument as in the proof of
Lemma \ref{lem-3} shows $g_2^{-1}g_1\in NT$, hence $g_2^{-1}g_1\in
B\cap NT=T$ and it follows that $g_1T=g_2T$. Now we show that
$\gamma$ is surjective. For $(x,gB)\in\varphi^{-1}(Y)$, we have
$x\in Y,\Ad(g^{-1})(x)\in\Lb$, hence there exists $b\in
B,x_0\in\Lt_0$ such that $\Ad(g^{-1})(x)=\Ad(b)(x_0)$. Then
$\gamma(x,gbT)=(x,gB)$.
\end{proof}

Since $\pi:\widetilde{Y}\rightarrow Y$ is a finite covering,
$\pi_!\bar{\mathbb{Q}}_{l\widetilde{Y}}$ is a well defined local
system on $Y$. Thus the intersection cohomology complex
$IC(\Lg,\pi_!\bar{\mathbb{Q}}_{l\widetilde{Y}})$ is well defined.
\begin{proposition}\label{p-2}
$\varphi_!\bar{\mathbb{Q}}_{lX}$ is canonically isomorphic to
$IC(\Lg,\pi_!\bar{\mathbb{Q}}_{l\widetilde{Y}})$. Moreover,
$\End(\varphi_!\bar{\mathbb{Q}}_{lX})=\End(\pi_!\bar{\mathbb{Q}}_{l\widetilde{Y}})=\bar{\bQ}_l[W]$.
\end{proposition}
\begin{proof}
Note that we have a commutative diagram
$$\CD
  \widetilde{Y} @>i\cdot\gamma>> X \\
  @V \pi VV @V \varphi VV  \\
  Y @>j>> \Lg,
\endCD$$ where $i:\varphi^{-1}(Y)\rightarrow X$ and $j$ are
inclusions. 
By base change theorem, we have
$\varphi_!\bar{\mathbb{Q}}_{lX}|_Y=\pi_!\bar{\mathbb{Q}}_{l\widetilde{Y}}$.
Since $\varphi$ is proper and $X$ is smooth of dimension equal to
$\dim Y$, we have that the Verdier dual
$\mathfrak{D}(\varphi_!\bar{\mathbb{Q}}_{lX})=\varphi_!(\mathfrak{D}\bar{\mathbb{Q}}_{lX})\cong\varphi_!\bar{\bQ}_{lX}[2\dim
Y]$. Hence by the definition of intersection cohomology complex, it
is enough to prove that
$$
\forall\
i>0,\dim\text{supp}\mathcal{H}^i(\varphi_!\bar{\mathbb{Q}}_{lX})<\dim
Y-i.
$$
For $x\in\Lg$, the stalk
$\cH^i_x(\varphi_!\bar{\mathbb{Q}}_{lX})=H^i_c(\varphi^{-1}(x),\bar{\mathbb{Q}}_{l})$.
Hence it is enough to show $\forall\
i>0,\dim\{x\in\Lg|H^i_c(\varphi^{-1}(x),\bar{\mathbb{Q}}_{l})\neq
0\}<\dim Y-i$. If $H^i_c(\varphi^{-1}(x),\bar{\mathbb{Q}}_{l})\neq
0$, then $i\leq 2\dim \varphi^{-1}(x)$. Hence it is enough to show
that $$\forall\ i>0,\dim\{x\in\Lg|\dim\varphi^{-1}(x)\geq i/2\}<\dim
Y-i.$$ Suppose this is not true for some $i$, then
$\dim\{x\in\Lg|\dim\varphi^{-1}(x)\geq i/2\}\geq\dim Y-i$. Let
$V=\{x\in\Lg|\dim\varphi^{-1}(x)\geq i/2\}$, it is closed in $\Lg$
but not equal to $\Lg$. Consider the map $p:Z\rightarrow \Lg,\
(x,B_1,B_2)\mapsto x$. We have $\dim p^{-1}(V)=\dim
V+2\dim\varphi^{-1}(x)\geq\dim V+i\geq\dim Y$ (for some $x\in V$).
 Thus by Lemma \ref{prop-dim}
(iii), $p^{-1}(V)$ contains some $Z_\mathcal{O},\mathcal{O}=G$-orbit
of $(B,nBn^{-1})$ in $\cB\times\cB$. If $x\in \Lt_0$, then
$(x,B,nBn^{-1})\in Z_{\mathcal{O}}$, hence $x$ belongs to the
projection of $p^{-1}(V)$ to $\Lg$ which has dimension $\dim V<\dim
Y$. But this projection is $G$-invariant hence contains all $Y$. We
get a contradiction.

Since $\pi$ is a principal $W$-bundle, we have
$\End(\pi_!\bar{\mathbb{Q}}_{l\widetilde{Y}})=\bar{\bQ}_l[W]$. It
follows that $\End(\varphi_!\bar{\mathbb{Q}}_{lX})=\bar{\bQ}_l[W]$.
\end{proof}

\subsection{}
In this subsection, we introduce some sheaves on the variety of
semisimple $G$-conjugacy classes in $\Lg$ similar to \cite{Lu1,Lu4}.

Let $\A$ be the set of closed $G$-conjugacy classes in $\Lg$. These
are precisely the semisimple classes in $\Lg$ (for a proof in the
group case see for example \cite{st}, and one can prove for the Lie
algebra case similarly). By geometric invariant theory, $\A$ has a
natural structure of affine variety and there is a well-defined
morphism $\sigma:\Lg\rightarrow\A$ such that $\sigma(x)$ is the
$G$-conjugacy class of $x_s$. There is a unique $\varsigma\in\A$
such that $\sigma^{-1}(\varsigma)=\{x\in\Lg|x \text{ nilpotent}\}$.

Recall that $Z=\{(x,B_1,B_2)\in
\mathfrak{g}\times\cB\times\cB|x\in\mathfrak{b}_1\cap\mathfrak{b}_2\}$.
Define $\tilde{\sigma}: Z\rightarrow\A$ by
$\tilde{\sigma}(x,B_1,B_2)=\sigma(x)$. For $a\in\A$, let
$Z^a=\tilde{\sigma}^{-1}(a)$.
\begin{lemma}
 We have $\dim Z^a\leq d_0$, where $d_0=\dim G-r$.
\end{lemma}
\begin{proof}
Define $m:Z^a\rightarrow\sigma^{-1}(a)$ by $(x,B_1,B_2)\mapsto x$.
Let $\rc\subset\sigma^{-1}(a)$ be a $G$-conjugacy class. Consider
$m:m^{-1}(\rc)\rightarrow\rc$. We have $\dim m^{-1}(\rc)\leq\dim
\rc+2(\dim G-r-\dim\rc)/2=\dim G-r$ (use Lemma \ref{prop-dim} (ii)).
Since $\sigma^{-1}(a)$ is a union of finitely many $G$-conjugacy
classes, it follows that $\dim Z^a\leq d_0$.
\end{proof}

Let $\mathcal{T}=\cH^{2d_0}\tilde{\sigma}_!\bar{\bQ}_{lZ}$. Recall
that we set $ Z_\cO=\{(x,B_1,B_2)\in Z|(B_1,B_2)\in \cO\}$, where
$\cO$ is an orbit of $G$ action on $\cB\times \cB$. Let
$\mathcal{T}^\cO=\cH^{2d_0}\sigma^0_!\bar{\bQ}_l$, where $\sigma^0:\
 Z_\cO\rightarrow\A$ is the restriction of $\tilde{\sigma}$ on
$Z_\cO$.

\begin{lemma}\label{lc-1}
We have $\cT^\cO\cong\bar{\sigma}_!\bar{\bQ}_l$, where
$\bar{\sigma}:\Lt\rightarrow\A$ is the restriction of $\sigma$.
\end{lemma}

\begin{proof}
The fiber of the natural projection $pr_{23}:Z_\cO\rightarrow \cO$
at $(B,nBn^{-1})\in \cO$ can be identified with $V=\Lb\cap n\Lb
n^{-1}$. Let $\cT'^\cO=\cH^{2d_0-2\dim \cO}\sigma'_!\bar{\bQ}_l$,
where $\sigma':V\rightarrow\A$ is $x\mapsto\sigma(x)$. Let
$\cT''^{\cO}=\cH^{2d_0+2\dim H}\sigma''_!\bar{\bQ}_l$, where
$H=B\cap nBn^{-1}$ and $\sigma'':G\times V\rightarrow\A$ is
$(g,x)\mapsto\sigma(x)$. Consider the composition $G\times
V\xrightarrow{pr_2}V\xrightarrow{\sigma'}\A$ (equal to $\sigma''$)
and the composition $G\times V\xrightarrow{p} H\backslash(G\times
V)=Z_\cO\xrightarrow{\sigma^0}\A$ (equal to $\sigma''$), we obtain
$$\cT''^{\cO}=\cH^{2d_0+2\dim
H}(\sigma'_!pr_{2!}\bar{\bQ}_l)=\cH^{2d_0+2\dim
H}(\sigma'_!\bar{\bQ}_l[-2\dim
G])=\cT'^{\cO},$$$$\cT''^{\cO}=\cH^{2d_0+2\dim
H}(\sigma^0_!p_{!}\bar{\bQ}_l)=\cH^{2d_0+2\dim
H}(\sigma^0_!\bar{\bQ}_l[-2\dim H])=\cT^{\cO}.$$ It follows that
$\cT^{\cO}=\cT'^{\cO}$. Now $V$ is fibred over $\Lt$ with fibers
isomorphic to $\Ln\cap n\Ln n^{-1}$. The map
$\sigma':V\rightarrow\A$ factors through
$\bar{\sigma}:\Lt\rightarrow \A$. Since $\Ln\cap n\Ln n^{-1}$ is an
affine space of dimension $d_0-\dim \cO$, we see that
$\cT^{\cO}=\cT'^{\cO}\cong\cH^0\bar{\sigma}_!\bar{\bQ}_l$. Since
$\bar{\sigma}$ is a finite covering (Lemma \ref{l-fc}), we have
$\cT^{\cO}\cong\bar{\sigma}_!\bar{\bQ}_l$.
\end{proof}
\begin{lemma}\label{l-fc}
The map $\bar{\sigma}:\Lt\rightarrow \A$ is a finite covering.
\end{lemma}
\begin{proof}
We show that for $x_1,x_2\in\Lt$, if $\sigma(x_1)=\sigma(x_2)$, then
there exists $w\in W$ such that $x_2=\Ad(w)x_1$. Since
$\sigma(x_1)=\sigma(x_2)$, there exists $g\in G$ such that
$x_1=\Ad(g)(x_2)$. It follows that $Z_G^0(x_1)=gZ_G^0(x_2)g^{-1}$.
We have $T\subset Z_G^0(x_1)$ and $gTg^{-1}\subset Z_G^0(x_1)$.
Hence there exists $h\in Z_G^0(x_1)$ such that $hTh^{-1}=gTg^{-1}$.
Let $n=g^{-1}h$, we have
$x_2=\Ad(g^{-1})x_1=\Ad(nh^{-1})x_1=\Ad(n)x_1$.
\end{proof}

Denote $\cT_\varsigma$ and $\cT^{\cO}_\varsigma$ the stalk  of $\cT$
and $\cT^{\cO}$ at $\varsigma$ respectively.
\begin{lemma}\label{p-1}
For $w\in W$, let $\cO_w$ be the $G$-orbit on $\cB\times\cB$ which
contains $(B,n_wBn_w^{-1})$. There is a canonical isomorphism
$\mathcal{T}_\varsigma\cong\bigoplus_{w\in
W}\mathcal{T}^{\cO_w}_\varsigma$.
\end{lemma}
\begin{proof}
We have
$\tilde{\sigma}^{-1}(\varsigma)=Z'=\{(x,B_1,B_2)\in\Lg\times\cB\times\cB|
x\in\Ln_1\cap\Ln_2\}$. We have a partition $Z'=\sqcup_{w\in
W}Z'_{\cO_w}$, where $Z'_{\cO_w}=\{(x,B_1,B_2)\in Z'|(B_1,B_2)\in
\cO_w\}$. Since $\dim Z'=d_0$, we have an isomorphism
\begin{equation*}
H^{2d_0}_c(Z',\bar{\bQ}_{l})= \bigoplus_{w\in
W}H^{2d_0}_c(Z'_{\cO_w},\bar{\bQ}_{l}),
\end{equation*}
which is $\mathcal{T}_\varsigma\cong\bigoplus_{w\in
W}\mathcal{T}_\varsigma^{\cO_w}$.
\end{proof}

Recall that we have
$\bar{\bQ}_l[W]=\End(\pi_!\bar{\bQ}_{l\widetilde{Y}})=\End(\varphi_!\bar{\bQ}_{lX})$.
In particular, $\varphi_!\bar{\bQ}_{lX}$ is naturally a $W$-module
and $\varphi_!\bar{\bQ}_{lX}\otimes\varphi_!\bar{\bQ}_{lX} $ is
naturally a $W$-module (with $W$ acting on the first factor). This
induces a $W$-module structure on
$\cH^{2d_0}\sigma_!(\varphi_!\bar{\bQ}_{lX}\otimes\varphi_!\bar{\bQ}_{lX})=\cT$.
Hence we obtain a $W$-module structure on the stalk $\cT_\varsigma$.

\begin{lemma}\label{l-multi}
Let $w\in W$. Multiplication by $w$ in the $W$-module structure of
$\cT_\varsigma=\bigoplus_{w'\in W}\cT_\varsigma^{\cO_{w'}}$ defines
for any $w'\in W$ an isomorphism
$\cT_\varsigma^{\cO_{w'}}\xrightarrow{\sim}\cT^{\cO_{ww'}}_\varsigma$.
\end{lemma}

\begin{proof}
We have an isomorphism
\begin{equation*}
f:Z'_{\cO_{w'}}\xrightarrow{\sim}Z'_{\cO_{ww'}},
(x,gBg^{-1},gn_{w'}Bn_{w'}^{-1}g^{-1})\mapsto(x,gn_w^{-1}Bn_wg^{-1},gn_{w'}Bn_{w'}^{-1}g^{-1}).
\end{equation*}
This induces an isomorphism
$$H^{2d_0}_c(Z'_{\cO_{w'}},\bar{\bQ}_{l})\xrightarrow{\sim}
H^{2d_0}_c(Z'_{\cO_{ww'}},\bar{\bQ}_{l})$$ which is just
multiplication by $w$.
\end{proof}

\subsection{}\label{ssec-3}
Let $\hat W$ be the set of simple modules (up to isomorphism) for
the Weyl group $W$ of $G$ (A description of $\hat W$ is given for
example in \cite{Lu3}). Given a semisimple object $M$ of some
abelian category such that $M$ is a $W$-module, we write
$M_\rho=\Hom_{\bar{\bQ}_l[W]}(\rho,M)$ for $\rho\in \Hat{W}$. We
have $M=\oplus_{\rho\in \Hat{W}}(\rho\otimes M_\rho)$ with $W$
acting on the $\rho$-factor and $M_\rho$ is in our abelian category.
In particular, we have $$
\pi_!\bar{\bQ}_{l\widetilde{Y}}=\bigoplus_{\rho\in
\Hat{W}}(\rho\otimes(\pi_!\bar{\bQ}_{l\widetilde{Y}})_{\rho}),$$
where $(\pi_!\bar{\bQ}_{l\widetilde{Y}})_{\rho}$ is an irreducible
local system on $Y$. We have $$
\varphi_!\bar{\bQ}_{lX}=\bigoplus_{\rho\in
\Hat{W}}(\rho\otimes(\varphi_!\bar{\bQ}_{lX})_{\rho}), $$ where
$(\varphi_!\bar{\bQ}_{lX})_{\rho}=IC(\bar{Y},(\pi_!\bar{\bQ}_{l\widetilde{Y}})_{\rho})$.
Moreover, for $a\in \A$, we have $\cT_a=\bigoplus_{\rho\in
\Hat{W}}(\rho\otimes(\cT_a)_{\rho})$. Set
$$\bar{Y}^\varsigma=\{x\in\bar{Y}|\sigma(x)=\varsigma\}, X^\varsigma=\varphi^{-1}(\bar{Y}^\varsigma)\subset X.$$ We have
$\bar{Y}^\varsigma=\{x\in\Lg|x \text{ nilpotent}\}$. Let
$\varphi^\varsigma:X^\varsigma\rightarrow\bar{Y}^\varsigma$ be the
restriction of $\varphi:X\rightarrow\Lg$.

\begin{lemma}\label{l-1}
$\mathrm{(i)}$ $X^\varsigma$ and $\bar{Y}^\varsigma$ are irreducible
varieties of dimension $d_0=\dim G-r$.

$\mathrm{(ii)}$ We have
$(\varphi_!\bar{\bQ}_{lX})|_{\bar{Y}^\varsigma}=\varphi^\varsigma_!\bar{\bQ}_{lX^\varsigma}$.
Moreover, $\varphi^\varsigma_!\bar{\bQ}_{lX^\varsigma}[d_0]$ is a
semisimple perverse sheaf on $\bar{Y}^\varsigma$.

$\mathrm{(iii)}$ We have
$(\varphi_!\bar{\bQ}_{lX})_\rho|_{\bar{Y}^\varsigma}\neq 0$ for any
$\rho\in \Hat{W}$.
\end{lemma}
\begin{proof}
(i) $\bar{Y}^\varsigma$, which is the nilpotent variety, is
well-known to be irreducible of dimension $\dim G-r$. We have
$X^\varsigma =\{(x,gB)\in\Lg\times\cB|\Ad(g^{-1})(x)\in\Ln\}$. By
projection to the second coordinate, we see that $\dim
X^\varsigma=\dim\Ln+\dim\cB=\dim G-r$. This proves (i).

The first assertion of (ii) follows by applying base change theorem
to the following commutative diagram$$\CD
  X^\varsigma @>i_1>> X \\
  @V \varphi^\varsigma VV @V \varphi VV  \\
  \bar{Y}^\varsigma @>i_2>> \Lg,
\endCD$$ where $i_1,i_2$ are inclusions.  Since $\varphi^\varsigma$
is proper, by similar argument as in the proof of Proposition
\ref{p-2}, to show that
$\varphi^\varsigma_!\bar{\bQ}_{lX^\varsigma}[d_0]$ is a perverse
sheaf, it suffices to show $$\forall\ i\geq
0,\dim\supp\cH^i(\varphi^\varsigma_!\bar{\bQ}_{lX^\varsigma})\leq
\dim\bar{Y}^\varsigma-i.$$ It is enough to show $\forall\ i\geq 0$,
$\dim\{x\in\bar{Y}^\varsigma|\dim(\varphi^\varsigma)^{-1}(x)\geq
i/2\}\leq\dim\bar{Y}^\varsigma-i$. If this is not true for some
$i\geq 0$, it would follow that the variety $\{(x,B_1,B_2)\in
\mathfrak{g}\times\cB\times\cB|x\in\mathfrak{n}_1\cap\mathfrak{n}_2\}$
has dimension greater than $\dim\bar{Y}^\varsigma=\dim G-r$, which
contradicts to Lemma \ref{prop-dim}. This proves that
$\varphi^\varsigma_!\bar{\bQ}_{lX^\varsigma}[d_0]$ is a perverse
sheaf. It is semisimple by the decomposition theorem\cite{BBD}. This
proves (ii).

Now we prove (iii). By Lemma \ref{lc-1}, we have
$\cT_\varsigma^{\cO_1}=H^0_c(\Lt\cap\sigma^{-1}(\varsigma),\bar{\bQ}_{l})\neq
0$. From Lemma \ref{l-multi}, we see that the $W$-module structure
defines an injective map
$\bar{\bQ}_l[W]\otimes\cT_\varsigma^{\cO_1}\rightarrow
\cT_\varsigma$. Since $\cT_\varsigma^{\cO_1}\neq 0$, we have
$(\bar{\bQ}_l[W]\otimes\cT_\varsigma^{\cO_1})_\rho\neq 0$ for any
$\rho\in \Hat{W}$, hence $(\cT_\varsigma)_\rho\neq 0$. We have
$\cT_\varsigma=H^{2d_0}_c(\bar{Y}^\varsigma,\varphi_!\bar{\bQ}_{lX}\otimes\varphi_!\bar{\bQ}_{lX})$,
hence
$$\bigoplus_{\rho\in \Hat{W}}\rho\otimes(\cT_\varsigma)_\rho=\bigoplus_{\rho\in \Hat{W}}\rho\otimes
H^{2d_0}_c(\bar{Y}^\varsigma,(\varphi_!\bar{\bQ}_{lX})_\rho\otimes\varphi_!\bar{\bQ}_{lX}).$$
This implies that
$(\cT_\varsigma)_{\rho}=H^{2d_0}_c(\bar{Y}^\varsigma,(\varphi_!\bar{\bQ}_{lX})_\rho\otimes\varphi_!\bar{\bQ}_{lX})$.
Thus it follows from $(\cT_\varsigma)_{\rho}\neq 0$ that
$(\varphi_!\bar{\bQ}_{lX})_\rho|_{\bar{Y}^\varsigma}\neq 0$ for any
$\rho\in \Hat{W}$.
\end{proof}

Let $\mathfrak{A}$ be the set of all pairs
$(\mathrm{c},\mathcal{F})$ where $\mathrm{c}$ is a nilpotent
$G$-conjugacy class in $\Lg$ and $\mathcal{F}$ is an irreducible
$G$-equivariant local system on $\mathrm{c}$ (up to isomorphism).

\begin{proposition}\label{mp-1}

$\mathrm{(i)}$ The restriction map
$\End_{\mathcal{D}(\bar{Y})}(\varphi_!\bar{\bQ}_{lX})
\rightarrow\End_{\mathcal{D}(\bar{Y}^\varsigma)}(\varphi^\varsigma_!\bar{\bQ}_{lX^\varsigma})$
is an isomorphism.

$\mathrm{(ii)}$ For any $\rho\in \Hat{W}$, there is a unique
$(\mathrm{c},\mathcal{F})\in\mathfrak{A}$ such that
$(\varphi_!\bar{\bQ}_{lX})_\rho|_{\bar{Y}^\varsigma}[d_0]$ is
$IC(\bar{\rc},\cF)[\dim\rc]$ regarded as a simple perverse sheaf on
$\bar{Y}^\varsigma$ (zero outside $\bar{\rc}$). Moreover,
$\rho\mapsto (\rc,\cF)$ is an injective map $\gamma:
\Hat{W}\rightarrow\mathfrak{A}$.
\end{proposition}
\begin{proof}
(i). Recall that we have $\varphi_!\bar{\bQ}_{lX}=\bigoplus_{\rho\in
\Hat{W}}\rho\otimes(\varphi_!\bar{\bQ}_{lX})_\rho$ where
$(\varphi_!\bar{\bQ}_{lX})_\rho[\dim\bar{Y}]$ are simple perverse
sheaves on $\bar{Y}$. Thus we have
$\varphi_!\bar{\bQ}_{lX}|_{\bar{Y}^\varsigma}=\varphi^\varsigma_!\bar{\bQ}_{lX^\varsigma}=\bigoplus_{\rho\in
\Hat{W}}\rho\otimes(\varphi_!\bar{\bQ}_{lX})_\rho|_{\bar{Y}^\varsigma}$
(we use Lemma \ref{l-1} (ii)). The restriction map
$\End_{\mathcal{D}(\bar{Y})}(\varphi_!\bar{\bQ}_{lX})\rightarrow
\End_{\mathcal{D}(\bar{Y}^\varsigma)}(\varphi^\varsigma_!\bar{\bQ}_{lX^\varsigma})$
is factorized as
\begin{eqnarray*}
&&\bigoplus_{\rho\in
\Hat{W}}\End_{\mathcal{D}(\bar{Y})}(\rho\otimes(\varphi_!\bar{\bQ}_{lX})_\rho)\xrightarrow{b}\bigoplus_{\rho\in
\Hat{W}}\End_{\mathcal{D}(\bar{Y}^\varsigma)}(\rho\otimes(\varphi_!\bar{\bQ}_{lX})_\rho|_{\bar{Y}^\varsigma})
\xrightarrow{c}\End_{\mathcal{D}(\bar{Y}^\varsigma)}(\varphi^\varsigma_!\bar{\bQ}_{lX^\varsigma})
\end{eqnarray*}
where $b=\oplus_\rho b_\rho$,
$b_\rho:\End(\rho)\otimes\End_{\mathcal{D}(\bar{Y})}((\varphi_!\bar{\bQ}_{lX})_\rho)
\rightarrow\End(\rho)\otimes\End_{\mathcal{D}(\bar{Y}^\varsigma)}((\varphi_!\bar{\bQ}_{lX})_\rho|_{\bar{Y}^\varsigma}).$
By Lemma \ref{l-1} (iii),
$(\varphi_!\bar{\bQ}_{lX})_\rho|_{\bar{Y}^\varsigma}\neq 0$, thus
$\End_{\mathcal{D}(\bar{Y})}((\varphi_!\bar{\bQ}_{lX})_\rho)=\bar{\bQ}_l
\subset\End_{\mathcal{D}(\bar{Y}^\varsigma)}((\varphi_!\bar{\bQ}_{lX})_\rho|_{\bar{Y}^\varsigma})$.
It follows that $b_\rho$ and thus $b$ is injective. Since $c$ is
also injective, the restriction map in injective. Hence it remains
to show that
$$\dim\End_{\mathcal{D}(\bar{Y}^\varsigma)}(\varphi^\varsigma_!\bar{\bQ}_{lX^\varsigma})=\dim\End_{\mathcal{D}(\bar{Y})}(\varphi_!\bar{\bQ}_{lX}).$$
For $A,A'$ two simple perverse sheaves on a variety $X$, we have
$H_c^0(X,A\otimes A')=0$ if and only if $A$ is not isomorphic to
$\mathfrak{D}(A')$ and $\dim H_c^0(X,A\otimes \mathfrak{D}(A))=1$
(see \cite{Lu2} section 7.4). We apply this to the semisimple
perverse sheaf $\varphi^\varsigma_!\bar{\bQ}_{lX^\varsigma}[d_0]$ on
$\bar{Y}^\varsigma$ and get
$\dim\End_{\mathcal{D}(\bar{Y}^\varsigma)}(\varphi^\varsigma_!\bar{\bQ}_{lX^\varsigma})=\dim
H^0_c(\bar{Y}^\varsigma,\varphi^\varsigma_!\bar{\bQ}_{lX^\varsigma}[d_0]\otimes\mathfrak{D}(\varphi^\varsigma_!\bar{\bQ}_{lX^\varsigma}[d_0]))$.
We have
\begin{eqnarray*}
&&\dim
H^0_c(\bar{Y}^\varsigma,\varphi^\varsigma_!\bar{\bQ}_{lX^\varsigma}[d_0]\otimes\mathfrak{D}(\varphi^\varsigma_!\bar{\bQ}_{lX^\varsigma}[d_0]))=\dim
H^0_c(\bar{Y}^\varsigma,\varphi^\varsigma_!\bar{\bQ}_{lX^\varsigma}[d_0]\otimes\varphi^\varsigma_!\bar{\bQ}_{lX^\varsigma}[d_0])\\&&
=\dim
H^{2d_0}_c(\bar{Y}^\varsigma,\varphi^\varsigma_!\bar{\bQ}_{lX^\varsigma}\otimes\varphi^\varsigma_!\bar{\bQ}_{lX^\varsigma})=\dim
H^{2d_0}_c(\bar{Y}^\varsigma,\varphi_!\bar{\bQ}_{lX}\otimes\varphi_!\bar{\bQ}_{lX})=\dim\cT_\varsigma=\sum_{w\in
W}\dim \cT_\varsigma^{\cO_w}.
\end{eqnarray*}
(The third equality follows from Lemma \ref{l-1} (ii) and the last
one follows from Lemma \ref{p-1}.)

We have
${}\cT_\varsigma^{\cO_w}=H^0_c(\bar{\sigma}^{-1}(\varsigma),\bar{\bQ}_l)$
(see Lemma \ref{lc-1}), hence $\dim {}\cT_\varsigma^{\cO_w}=1$ and
$$\sum_{w\in W}\dim \cT_\varsigma^{\cO_w}=|W|=\dim\End_{\mathcal{D}(\bar{Y})}(\varphi_!\bar{\bQ}_{lX}).$$
Thus (i) is proved.

From the proof of (i) we see that both $b$ and $c$ are isomorphisms.
It follows that the perverse sheaf
$(\varphi_!\bar{\bQ}_{lX})_\rho|_{\bar{Y}^\varsigma}[d_0]$ on
$\bar{Y}^\varsigma$ is simple and that for $\rho,\rho'\in \Hat{W}$,
we have
$(\varphi_!\bar{\bQ}_{lX})_\rho|_{\bar{Y}^\varsigma}[d_0]\cong(\varphi_!\bar{\bQ}_{lX})_{\rho'}|_{\bar{Y}^\varsigma}[d_0]$
if and only if $\rho=\rho'$. Since the simple perverse sheaf
$(\varphi_!\bar{\bQ}_{lX})_\rho|_{\bar{Y}^\varsigma}[d_0]$ is
$G$-equivariant and $\bar{Y}^\varsigma$ consists of finitely many
nilpotent $G$-conjugacy classes,
$(\varphi_!\bar{\bQ}_{lX})_\rho|_{\bar{Y}^\varsigma}[d_0]$ must be
as in (ii).
\end{proof}

\subsection{}

In this subsection let $G=SO_N(\tk)$ (resp. $Sp_{2n}(\tk)$) and
$\Lg=\Lo_N(\tk)$ (resp. $\mathfrak{sp}_{2n}(\tk)$) be the Lie
algebra of $G$. Let $G_{ad}$ be an adjoint group over $\tk$ of the
same type as $G$ and $\Lg_{ad}$ be the Lie algebra of $G_{ad}$. For
$q$ a power of 2, let $G(\tF_q)$, $\mathfrak{g}({\tF}_q)$ be the
fixed points of a split Frobenius map $\mathfrak{F}_q$ relative to
$\tF_q$ on $G$, $\Lg$. Let $G_{ad}(\tF_q)$, $\Lg_{ad}(\tF_q)$ be
defined like $G(\tF_q)$, $\Lg(\tF_q)$. Let $\mathfrak{A}$ be the set
of all pairs $(\mathrm{c},\mathcal{F})$ where $\mathrm{c}$ is a
nilpotent $G$-orbit in $\Lg$ and $\mathcal{F}$ is an irreducible
$G$-equivariant local system on $\mathrm{c}$ (up to isomorphism).
Let $\mathfrak{A}_{ad}$ be defined for $G_{ad}$ as in the
introduction. We show that the number of elements in
$\mathfrak{A}_{ad}$ is equal to the number of elements in
$\mathfrak{A}$.

We first show that the number of elements in $\mathfrak{A}$ is equal
to the number of nilpotent $G(\tF_q)$-orbits in
$\mathfrak{g}({\tF}_q)$ (for $q$ large). To see this we can assume
$\tk=\bar{\tF}_2$.   Pick representatives $x_1,\ldots,x_M$ for the
nilpotent $G$-orbits in $\Lg$. If $q$ is large enough, the Frobenius
map $\mathfrak{F}_q$ keeps $x_i$ fixed and acts trivially on
$Z_{G}(x_i)/Z_{G}^0(x_i)$. Then the number of $G(\tF_q)$-orbits in
the $G$-orbit of $x_i$ is equal to the number of irreducible
representations of $Z_{G}(x_i)/Z_{G}^0(x_i)$ hence to the number of
$G$-equivariant irreducible local systems on the $G$-orbit of $x_i$.
Similarly, the number of elements in $\mathfrak{A}_{ad}$ is equal to
the number of nilpotent $G_{ad}(\tF_q)$-orbits in
$\mathfrak{g}_{ad}({\tF}_q)$.

On the other hand, the number of nilpotent $G(\tF_q)$-orbits in
$\mathfrak{g}({\tF}_q)$ is equal to the number of nilpotent
$G_{ad}(\tF_q)$-orbits in $\mathfrak{g}_{ad}({\tF}_q)$. In fact, we
have a morphism $G\rightarrow G_{ad}$ which is an isomorphism of
abstract groups and an obvious bijective morphism
$\mathcal{U}\rightarrow \mathcal{U}_{ad}$ between the nilpotent
variety $\mathcal{U}$ of $ \Lg$ and the nilpotent variety
$\mathcal{U}_{ad}$ of $\Lg_{ad}$. Thus the nilpotent orbits in $\Lg$
and $\Lg_{ad}$ are in bijection and the corresponding
 component groups of centralizers are isomorphic. It follows that
 $|\mathfrak{A}|=|\mathfrak{A}_{ad}|$.

\begin{corollary}\label{coro-1}
$|\mathfrak{A}|=|\mathfrak{A}_{ad}|=|\hat{W}|$.
\end{corollary}
\begin{proof}
Assume $G=SO_N(\tk)$. From Proposition \ref{prop-num} (i), Corollary
\ref{coro} and the above argument we see that
$|\mathfrak{A}|=|\mathfrak{A}_{ad}|\leq |\hat{W}|$. On the other
hand, by Proposition \ref{mp-1} (ii),we have
$|\mathfrak{A}_{ad}|\geq|\hat{W}|$.

Assume $G=Sp_{2n}(\tk)$. It is known in \cite{Spal} that the number
of nilpotent $G(\tF_q)$-orbits in $\Lg(\tF_q)$ is equal to
$|\hat{W}|$. The assertion follows from the above argument.
\end{proof}

\begin{theorem}\label{coro-3}
The map $\gamma$ in Proposition \ref{mp-1} $\mathrm{(ii)}$ is a
bijection.
\end{theorem}

\begin{corollary}\label{coro-2}
Proposition \ref{prop-2}, Corollary \ref{coro-n}, Proposition
\ref{prop-num}, Corollary \ref{coro} hold with all "at most"
removed.
\end{corollary}
\begin{proof}
For $q$ large enough, this follows from Corollary \ref{coro-1}. Now
let $q$ be an arbitrary power of $2$. The Frobenius map
$\mathfrak{F}_q$ acts trivially on $W$. Since we have the Springer
correspondence map $\gamma$ in Proposition \ref{mp-1}
$\mathrm{(ii)}$, for each pair $(\rc,\cF)\in\mathfrak{A}_{ad}$,
$\rc$ is stable under the Frobenius map $\mathfrak{F}_q$ and we have
$\mathfrak{F}_q^{-1}(\cF)\cong\cF$. Pick a rational point $x$ in
$\rc$. The $G_{ad}$-equivariant local systems on $\rc$ are in 1-1
correspondence with the isomorphism classes of the irreducible
representations of $Z_{G_{ad}}(x)/Z_{G_{ad}}^0(x)$. Since
$Z_{G_{ad}}(x)/Z_{G_{ad}}^0(x)$ is abelian (see Proposition
\ref{prop-cen}) and the Frobenius map $\mathfrak{F}_q$ acts
trivially on the irreducible representations of
$Z_{G_{ad}}(x)/Z_{G_{ad}}^0(x)$, $\mathfrak{F}_q$ acts trivially on
$Z_{G_{ad}}(x)/Z_{G_{ad}}^0(x)$. Thus it follows that the number of
nilpotent $G_{ad}(\tF_q)$-orbits in $\mathfrak{g}_{ad}({\tF}_q)$ is
independent of $q$ hence it is equal to
$|\mathfrak{A}_{ad}|=|\hat{W}|$.
\end{proof}

A corollary of Theorem \ref{coro-3} is that in this case there are
no cuspidal local systems similarly defined as in \cite{Lu1}. This
result does not extend to exceptional Lie algebras. (In type $F_4$,
characteristic 2, the results of \cite{Spa2} suggest that a cuspidal
local system exists on a nilpotent class.)
\section{component groups of centralizers}

In this section we describe the component groups of centralizers
$Z_G(x)/Z_G^0(x)$, where $G=SO_N(\tk)$ and $x\in\Lo_N(\tk)$ is
nilpotent.
\begin{proposition}\label{prop-cen}
$\mathrm{(i)}$ Assume $V$ is defective. Let $x\in \Lo(V)$ be a
nilpotent element corresponding to the form module
$(\lambda_1)_{\chi(\lambda_1)}^{m_1}(\lambda_2)_{\chi(\lambda_2)}^{m_2}
\cdots(\lambda_s)_{\chi(\lambda_s)}^{m_s}$.  Denote by $n_1$ the
cardinality of $\{1\leq i\leq
s-1|\chi(\lambda_i)+\chi(\lambda_{i+1})\leq \lambda_i,\
\chi(\lambda_i)\neq\lambda_i/2\}$. We have
$Z_{O(V)}(x)/Z_{O(V)}^0(x)=(\mathbb{Z}_2)^{n_1}$.

$\mathrm{(ii)}$ Assume $V$ is non-defective. Let $x\in \Lo(V)$ be a
nilpotent element corresponding to the form module
$(\lambda_1)_{\chi(\lambda_1)}^{m_1}(\lambda_2)_{\chi(\lambda_2)}^{m_2}
\cdots(\lambda_s)_{\chi(\lambda_s)}^{m_s}$.  Denote by $n_2$ the
cardinality of $\{1\leq i\leq
s|\chi(\lambda_i)+\chi(\lambda_{i+1})\leq \lambda_i,\
\chi(\lambda_i)\neq\lambda_i/2\}$ (here define
$\chi(\lambda_{s+1})=0$). Assume $n_2\geq 1$. Then
$Z_{SO(V)}(x)/Z_{SO(V)}^0(x)=(\mathbb{Z}_2)^{n_2-1}$.
\end{proposition}
\begin{proof}
$\mathrm{(i)}$ We write $Z=Z_{O(V)}(x)$ and $Z^0=Z_{O(V)}^0(x)$ for
simplicity. We can assume that $q$ is large enough. (i) is proved in
two steps.

Step 1: we show that $Z/Z^0$ is an abelian group of order $2^{n_1}$.

The group $Z/Z^0$ has $2^{n_1}$ conjugacy classes, since the
$G$-orbit of $x$ splits into $2^{n_1}$ $G(\tF_q)$-orbits in
$\Lg(\tF_q)$ (Corollary \ref{coro-2}). We show $|Z/Z^0|=2^{n_1}$ by
showing that
$$(*)\quad\quad|Z(\tF_q)|=2^{n_1}q^{\dim(Z)}+\text{ lower terms in }q.$$
We prove $(*)$ by induction on $n_1$. Suppose $n_1=0$, then $Z/Z^0$
has only one conjugacy class. It follows that $Z/Z^0=\{1\}$ and
$(*)$ holds for $n_1=0$. Suppose $n_1\geq 1$. Let $t$ be the minimal
integer such that $\chi(\lambda_t)\neq \lambda_t/2$ and $
\chi(\lambda_t)+\chi(\lambda_{t+1})\leq \lambda_t$. Let
$V_1=(\lambda_1)^{m_1}_{\chi(\lambda_1)}(\lambda_2)^{m_2}_{\chi(\lambda_2)}\cdots(\lambda_t)^{m_t}_{\chi(\lambda_t)}$
and
$V_2=(\lambda_{t+1})^{m_{t+1}}_{\chi(\lambda_{t+1})}\cdots(\lambda_s)^{m_s}_{\chi(\lambda_s)}$.
Then $V_1$ is non-defective. Write $Z_i=Z_{O(V_i)}(V_i)$ and
$Z_i^0=Z_{O(V_i)}^0(V_i)$, $i=1,2$. We have
$|Z_{SO(V_1)}(V_1)/Z_{SO(V_1)}^0(V_1)|=1\Rightarrow |Z_1/Z_1^0|=2$.
It follows that $|Z_1(\tF_q)|=2q^{\dim Z_1}+$ lower terms in $q$. We
show that $|Z(\tF_q)|=|Z_1(\tF_q)|\cdot|Z_2(\tF_q)|\cdot
q^{\dim\Hom_A(V_1,V_2)}$. Then the assertion $(*)$ follows from
induction hypothesis since we have $\dim Z_1+\dim Z_2+\dim
\Hom_A(V_1,V_2)=\dim Z$.

 Consider $V_1$
as an element in the Grassmannian variety $Gr(V,r)$ of dimension
$r=\sum_{j=1}^tm_j\lambda_j$. Then $C(V)=\{g\in GL(V)|gx=xg\}$ acts
on $Gr(V,r)$. We have $C(V)(V_1\oplus V_2)=C(V)V_1\oplus
C(V)V_2\cong V_1\oplus V_2$. By our choice of $V_1$ and $V_2$, it
follows that $C(V)V_1\cong V_1$ and $C(V)V_2\cong V_2$. Thus the
orbit of $V_1$ under $C(V)$ coincides with the orbit of $V_1$ under
the action of $Z$. It is easy to verify that this orbit consists of
$q^{\dim \Hom_A(V_1,V_2)}$ elements (using $C(V)$ action). Since the
stabilizer of $V_1$ in $Z$ is the product of $Z_1$ and $Z_2$, we get
$|Z(\tF_q)|=|Z_1(\tF_q)|\cdot|Z_2(\tF_q)|\cdot
q^{\dim\Hom_A(V_1,V_2)}.$

Step 2: we show that there is a subgroup
$(\mathbb{Z}_2)^{n_1}\subset Z/Z^0$. Thus $Z/Z^0$ has to be
$(\mathbb{Z}_2)^{n_1}$. Let $1\leq i_1,\ldots,i_{n_1}\leq s-1$ be
such that $\chi(\lambda_{i_j})>\lambda_{i_j}/2$ and
$\chi(\lambda_{i_j})+\chi(\lambda_{i_{j}+1})\leq \lambda_{i_j}$,
$j=1,\ldots,n_1$. Let
$V_j=(\lambda_{i_{j-1}+1})_{\chi(\lambda_{i_{j-1}+1})}^{m_{i_{j-1}+1}}\cdots
(\lambda_{i_{j}})_{\chi(\lambda_{i_{j}})}^{m_{i_{j}}}$,
$j=1,\ldots,n_1+1$, where $i_0=0$, $i_{n_1+1}=s$. Then $V=V_1\oplus
V_2\oplus\cdots\oplus V_{n_1+1}$, where $V_i$, $i=1,\ldots, n_1$ are
non-defective and $V_{n_1+1}$ is defective and non-degenerate. We
have $Z_{O(V_i)}(V_i)/Z_{O(V_i)}^0(V_i)=\mathbb{Z}_2$,
$i=1,\ldots,n_1$, and
$Z_{O(V_{n_1+1})}(V_{n_1+1})/Z_{O(V_{n_1+1})}^0(V_{n_1+1})=\{1\}$.
Take $g_i\in Z_{O(V_i)}(V_i)$ such that $g_iZ_{O(V_i)}^0(V_i)$
generates $Z_{O(V_i)}(V_i)/Z_{O(V_i)}^0(V_i)$, $i=1,\ldots,n_1$. We
know each $V_i$, $i=1,\ldots,n_1$, has two isomorphism classes
$V_i^0,V_i^\delta$ over $\tF_q$. Let $m_i^0$ and $m_i^\delta$ be two
elements in $\Lo(V_i)(\tF_q)$ corresponding to $V_i^0$ and
$V_i^\delta$ respectively. We can assume $x=m_1^0\oplus
m_2^0\oplus\cdots\oplus m_{n_1}^0\oplus m_{n_1+1}$.  Let
$\tilde{g_i}=Id\oplus\cdots\oplus g_i\oplus\cdots\oplus Id$,
$i=1,\ldots,n_1$. Then we have $\tilde{g_i}\in Z$ and
$\tilde{g_i}\notin Z^0$, since $V_1^0\oplus\cdots\oplus
V_i^0\oplus\cdots\oplus V_{n_1}^0\oplus V_{n_1+1}\ncong
V_1^0\oplus\cdots\oplus V_i^\delta\oplus\cdots\oplus V_{n_1}^0\oplus
V_{n_1+1}$(Corollary \ref{coro-2}). We also have that the images of
$\tilde{g}_{i_1}\cdots\tilde{g}_{i_p}$'s, $1\leq i_1<\cdots<i_p\leq
n_1$, $p=1,\ldots,n_1$, in $Z/Z^0$ are not equal to each other.
Moreover $\tilde{g}_i^2\in Z^0$. Thus the $\tilde{g}_iZ^0$'s
generate a subgroup $(\mathbb{Z}_2)^{n_1}$ in $Z/Z^0$.

$\mathrm{(ii)}$
 Let us write $Z=Z_{SO(V)}(x)$ and $Z^0=Z_{SO(V)}^0(x)$ for simplicity. Assume $n_2\geq 1$. We
know that the group $Z/Z^0$ has $2^{n_2-1}$ conjugacy classes since
the $G$-orbit of $x$ splits into $2^{n_2-1}$ $G(\tF_q)$-orbits in
$\Lg(\tF_q)$ (we can assume $q$ large enough). The same argument as
in $\mathrm{(i)}$ shows that $|Z(\tF_q)|=2^{n_2-1}q^{\dim(Z)}+$lower
terms. Then it follows that $Z/Z^0$ is an abelian group of order
$2^{n_2-1}$. It is enough to show that there is a subgroup
$(\mathbb{Z}_2)^{n_2-1}\subset Z/Z^0$. Let $1\leq
i_1,\ldots,i_{n_2}\leq s$ be such that
$\chi(\lambda_{i_j})>\lambda_{i_j}/2$ and
$\chi(\lambda_{i_j})+\chi(\lambda_{i_{j}+1})\leq \lambda_{i_j}$,
$j=1,\ldots,n_2$.

Case 1: $\chi(\lambda_s)=\lambda_s/2$. Then $i_{n_2}<s$. Let
$V_j=(\lambda_{i_{j-1}+1})_{\chi(\lambda_{i_{j-1}+1})}^{m_{i_{j-1}+1}}\cdots
(\lambda_{i_{j}})_{\chi(\lambda_{i_{j}})}^{m_{i_{j}}}$,
$j=1,\ldots,n_2+1$, where $i_0=0$, $i_{n_2+1}=s$. Then $V=V_1\oplus
\cdots\oplus V_{n_2+1}$. We have
$Z_{O(V_i)}(V_i)/Z_{O(V_i)}^0(V_i)=\mathbb{Z}_2$, $i=1,\ldots,n_2$
and
$Z_{O(V_{n_2+1})}(V_{n_2+1})/Z_{O(V_{n_2+1})}^0(V_{n_2+1})=\{1\}$.
Take $g_i\in Z_{O(V_i)}(V_i)$ such that $g_iZ_{O(V_i)}^0(V_i)$
generates $Z_{O(V_i)}(V_i)/Z_{O(V_i)}^0(V_i)$, $i=1,\ldots,n_2$. Let
$\tilde{g_i}=g_1\oplus Id\oplus\cdots\oplus g_i\oplus\cdots\oplus
Id$, $i=2,\ldots,n_2$.  We have $\tilde{g_i}\in Z$ and
$\tilde{g_i}\notin Z^0$. We also have the images of
$\tilde{g}_{i_1}\cdots\tilde{g}_{i_p}$'s, $2\leq i_1<\cdots<i_p\leq
n_2$, $p=1,\ldots,n_2-1$, in $Z/Z^0$ are not equal to each other.
Moreover $\tilde{g}_i^2\in Z^0$. Hence the $\tilde{g}_iZ^0$'s
generate a subgroup $(\mathbb{Z}_2)^{n_2-1}$ in $Z/Z^0$.

Case2: $\chi(\lambda_s)>\lambda_s/2$. Then $i_{n_2}=s$. Let
$V_j=(\lambda_{i_{j-1}+1})_{\chi(\lambda_{i_{j-1}+1})}^{m_{i_{j-1}+1}}\cdots
(\lambda_{i_{j}})_{\chi(\lambda_{i_{j}})}^{m_{i_{j}}}$,
$j=1,\ldots,n_2$, where $i_0=0$. Then $V=V_1\oplus \cdots\oplus
V_{n_2}$. We have $Z_{O(V_i)}(V_i)/Z_{O(V_i)}^0(V_i)=\mathbb{Z}_2$,
$i=1,\ldots,n_2$. Take $g_i\in Z_{O(V_i)}(V_i)$ such that
$g_iZ_{O(V_i)}^0(V_i)$ generates
$Z_{O(V_i)}(V_i)/Z_{O(V_i)}^0(V_i)$, $i=1,\ldots,n_2$. Let
$\tilde{g_i}=g_1\oplus Id\oplus\cdots\oplus g_i\oplus\cdots\oplus
Id$, $i=2,\cdots,n_2$. The $\tilde{g}_iZ^0$'s generate a subgroup
$(\mathbb{Z}_2)^{n_2-1}$ in $Z/Z^0$.
\end{proof}

\section{Complements}In this section we give an example in type $B_4$. We also make some
comments on the explicit Springer correspondence.

\subsection{ Example} We list the unipotent classes in $SO(9)$ and
nilpotent classes in $\mathfrak{so}(9)$. We use the notation in
\cite{Hes}, see also section 1.

\vspace{.1in}

Unipotent classes in characteristic 2:
\begin{eqnarray*}&&\underline{8_51_1},\ \underline{6_42_21_1},\
\underline{6_41^3_1},\ \underline{4^2_31_1},\ 4^2_21_1,\
{4_32_2^21_1},\ \underline{4_32^2_11_1},\ 4_32_21^3_1,\
\underline{4_31^5_1},\ 3^2_22_21_1,\ \underline{3^2_21^3_1},\\&&
2^4_21_1,\ 2^4_11_1,\ 2^3_21^3_1,\ 2^2_21^5_1,\ 2^2_11^5_1,\
2_21^7_1,\ 1^9_1.
\end{eqnarray*}

 Nilpotent classes in characteristic 2: \begin{eqnarray*}&&5_54_4,\
4^2_41_1,\ \underline{4^2_31_1},\ 4^2_21_1,\ 4_43_31_1^2,\
3^2_31^3_1,\ \underline{3^2_21^3_1},\ 3^2_32_21_1,\ 3^2_22_21_1,\
3_32^3_2,\ 3_32_21^4_1,\\&& 2^4_21_1,\ 2^4_11_1,\ 2^3_21^3_1,\
2^2_21^5_1,\ 2^2_11^5_1,\ 2_21^7_1,\ 1^9_1.
\end{eqnarray*}

Unipotent/nilpotent classes in characteristic not 2:
\begin{eqnarray*}&& 9_5,\ \underline{7_4^11_1^2},\ \underline{\underline{5_3^13_2^11_1^1}},\
5_3^12_1^2,\ \underline{5_3^11_1^4},\ 4_2^21_1^1,\ 3_2^3,\
\underline{3_2^21_1^3},\ \underline{3_2^12_1^21_1^2},\
\underline{3_2^11_1^6},\ 2_1^41_1^1,\ 2_1^21_1^5,\
1_1^9.\end{eqnarray*}
In each case the component group of the centralizer is trivial
except that
 the component groups for the underlined ones are
$\mathbb{Z}/2\mathbb{Z}$ and for the double underlined one is
$(\mathbb{Z}/2\mathbb{Z})^2$. In this case the number of unipotent
classes and nilpotent classes over an algebraically closed field of
characteristic 2 happen to be the same, but this is not true in
higher ranks. Note that the component groups are already quite
different.

\subsection{}In \cite{Spal}, Spaltenstein shows that the explicit Springer
correspondence for the pairs $(\rc,\bar{\bQ}_l)$ is given by the
pairs of partitions associated by him to each orbit, assuming the
theory of Springer representation is valid. We will show elsewhere
that the explicit correspondence for the pairs $(\rc,\cF)$ ($\cF$
nontrivial) is given by the pairs of partitions appearing in the
proof of Proposition \ref{prop-num}.

The explicit generalized Springer correspondence for the unipotent
case in characteristic 2 is described by Lusztig and Spaltenstein in
\cite {LS} for classical groups. The correspondence is more
complicated than the nilpotent case because of the existence of
cuspidal local systems. In the unipotent case various Weyl groups
for Levi subgroups are needed and in our case we only need the Weyl
group for $G$.

 \vskip 10pt {\noindent\bf\large Acknowledgement} \vskip 5pt I
would like to thank Professor George Lusztig for his guidance,
encouragement and help during this research. I am also grateful to
the referee for the valuable suggestions.

\end{document}